\documentclass[12pt]{article}

\newenvironment{wileykeywords}{\textsf{Keywords:}\hspace{\stretch{1}}}{\hspace{\stretch{1}}\rule{1ex}{1ex}}

\usepackage{amsmath,amssymb}
\usepackage{graphicx}
\usepackage{color}
\usepackage{dcolumn}
\usepackage{bm}
\usepackage{url}
\usepackage[numbers,super,comma,sort&compress]{natbib}
\usepackage{multirow}
\usepackage{multicol}
\usepackage{hyperref}

\definecolor{background-color}{gray}{0.98}

\newcommand{\dx}{\ensuremath{\,\mathrm{d}\mathbf{x}}}
\newcommand{\dy}{\ensuremath{\,\mathrm{d}\mathbf{y}}}
\newcommand{\uint}{\ensuremath{\phi_\text{int}}}
\newcommand{\uext}{\ensuremath{\phi_\text{ext}}}
\newcommand{\phiint}{\ensuremath{\phi_\text{int}}}
\newcommand{\phiext}{\ensuremath{\phi_\text{ext}}}
\newcommand{\epsilonint}{\ensuremath{\epsilon_\text{int}}}
\newcommand{\epsilonext}{\ensuremath{\epsilon_\text{ext}}}
\newcommand{\normal}{\ensuremath{\hat{\mathbf{n}}}}

\newcommand{\traceDe}{\ensuremath{\gamma_D^+}}
\newcommand{\traceDi}{\ensuremath{\gamma_D^-}}
\newcommand{\traceNe}{\ensuremath{\gamma_N^+}}
\newcommand{\traceNi}{\ensuremath{\gamma_N^-}}
\newcommand{\traceDei}{\ensuremath{\gamma_D^\pm}}
\newcommand{\traceNei}{\ensuremath{\gamma_N^\pm}}
\newcommand{\SLP}{\ensuremath{\mathcal{V}}}
\newcommand{\DLP}{\ensuremath{\mathcal{K}}}
\newcommand{\SLPint}{\ensuremath{\mathcal{V}_\mathrm{int}}}
\newcommand{\DLPint}{\ensuremath{\mathcal{K}_\mathrm{int}}}
\newcommand{\SLPext}{\ensuremath{\mathcal{V}_\mathrm{ext}}}
\newcommand{\DLPext}{\ensuremath{\mathcal{K}_\mathrm{ext}}}
\newcommand{\SLPei}{\ensuremath{\mathcal{V}_\mathrm{ext,int}}}
\newcommand{\DLPei}{\ensuremath{\mathcal{K}_\mathrm{ext,int}}}
\newcommand{\SL}{\ensuremath{V}}
\newcommand{\DL}{\ensuremath{K}}
\newcommand{\AD}{\ensuremath{T}}
\newcommand{\HS}{\ensuremath{D}}
\newcommand{\SLint}{\ensuremath{V_\mathrm{int}}}
\newcommand{\DLint}{\ensuremath{K_\mathrm{int}}}
\newcommand{\ADint}{\ensuremath{T_\mathrm{int}}}
\newcommand{\HSint}{\ensuremath{D_\mathrm{int}}}
\newcommand{\SLext}{\ensuremath{V_\mathrm{ext}}}
\newcommand{\DLext}{\ensuremath{K_\mathrm{ext}}}
\newcommand{\ADext}{\ensuremath{T_\mathrm{ext}}}
\newcommand{\HSext}{\ensuremath{D_\mathrm{ext}}}
\newcommand{\SLei}{\ensuremath{V_\mathrm{ext,int}}}
\newcommand{\DLei}{\ensuremath{K_\mathrm{ext,int}}}
\newcommand{\ADei}{\ensuremath{T_\mathrm{ext,int}}}
\newcommand{\HSei}{\ensuremath{D_\mathrm{ext,int}}}
\newcommand{\ID}{\ensuremath{I}}

\title{Towards optimal boundary integral formulations of the Poisson-Boltzmann equation for molecular electrostatics}
\author{Stefan D. Search\thanks{Department of Mechanical Engineering, Universidad T\'ecnica Federico Santa Mar\'ia}, 
Christopher D. Cooper\footnotemark[1]~\thanks{Centro Cient\'ifico Tecnol\'ogico de Valpara\'iso, Universidad T\'ecnica Federico Santa Mar\'ia},
Elwin van 't Wout\thanks{Institute for Mathematical and Computational Engineering, School of Engineering and Faculty of Mathematics, Pontificia Universidad Cat\'olica de Chile}}

\begin{document}

\maketitle

\begin{abstract}
The Poisson-Boltzmann equation offers an efficient way to study electrostatics in molecular settings. Its numerical solution with the boundary element method is widely used, as the complicated molecular surface is accurately represented by the mesh, and the point charges are  accounted for explicitly. In fact, there are several well-known boundary integral formulations available in the literature. This work presents a generalized expression of the boundary integral representation of the implicit solvent model, giving rise to new forms to compute the electrostatic potential. Moreover, it proposes a strategy to build efficient preconditioners for any of the resulting systems, improving the convergence of the linear solver. We perform systematic benchmarking of a set of formulations and preconditioners, focusing on the time to solution, matrix conditioning, and eigenvalue spectrum. We see that the eigenvalue clustering is a good indicator of the matrix conditioning, and show that they can be easily manipulated by scaling the preconditioner. Our results suggest that the optimal choice is problem-size dependent, where a simpler direct formulation is the fastest for small molecules, but more involved second-kind equations are better for larger problems. We also present a fast Calderón preconditioner for first-kind formulations, which shows promising behavior for future analysis. This work sets the basis towards choosing the most convenient boundary integral formulation of the Poisson-Boltzmann equation for a given problem.   

\end{abstract}

\begin{wileykeywords}
Boundary element method, Preconditioning, Poisson-Boltzmann, Implicit solvent model, Electrostatics.
\end{wileykeywords}

\clearpage


\begin{figure}[h]
\centering
\colorbox{background-color}{
\fbox{
\begin{minipage}{1.0\textwidth}
\includegraphics[width=50mm,height=50mm]{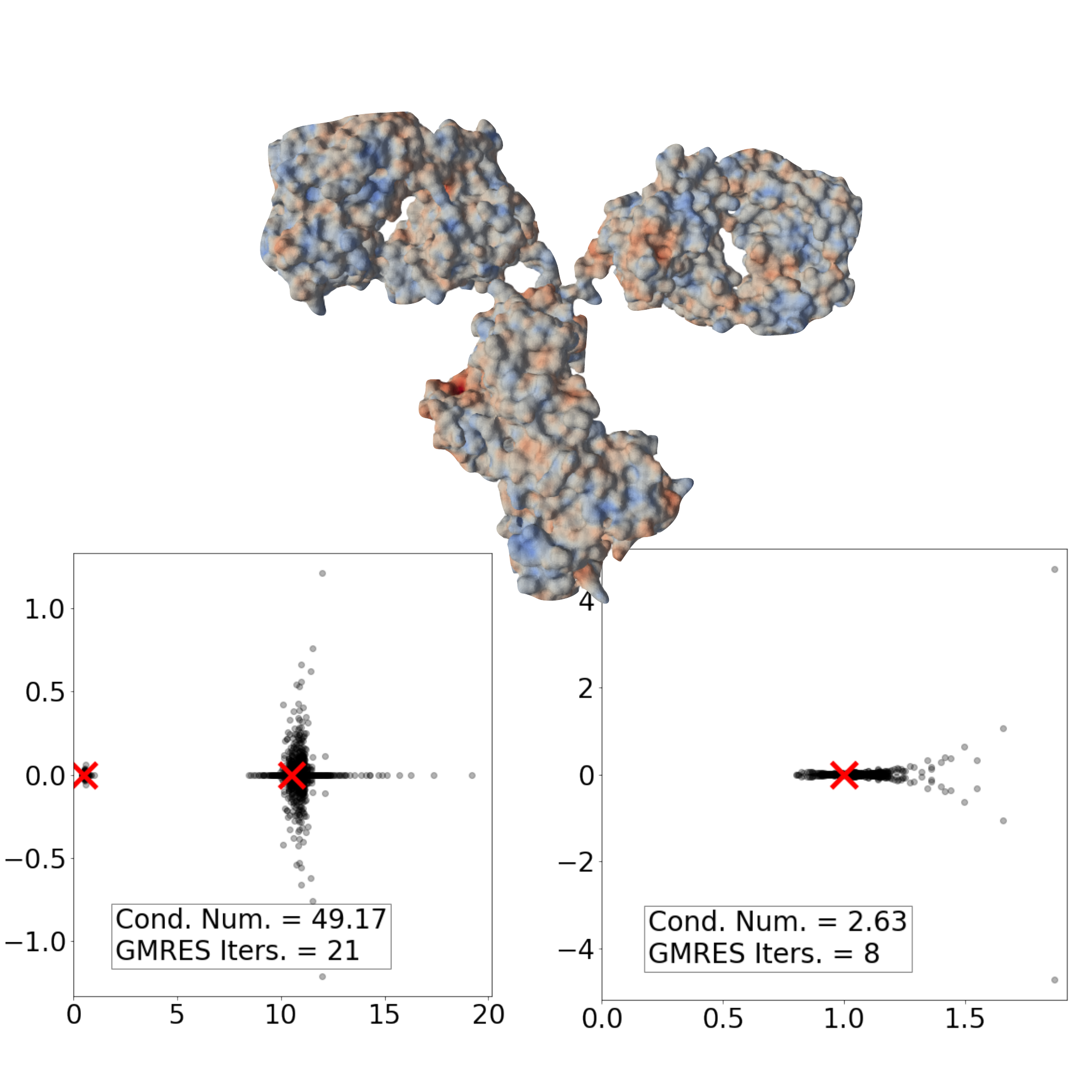} 
\\
The boundary element method is an efficient numerical technique to solve the Poisson-Boltzmann equation in molecular settings. This article presents a generalized form of the boundary integral equation, that leads to a family of expressions and preconditioners that are used here for the first time. Extensive benchmark tests show that choosing the right formulation has a huge impact in the time to solution, and that the optimal formulation depends on the molecule size.
\end{minipage}
}}
\end{figure}

  \makeatletter
  \renewcommand\@biblabel[1]{#1.}
  \makeatother

\bibliographystyle{apsrev}

\renewcommand{\baselinestretch}{1.5}
\normalsize

\clearpage

\section*{\sffamily \Large INTRODUCTION} 

\begin{figure}[!ht]
    \centering
    \includegraphics[width=0.3\textwidth,keepaspectratio=true]{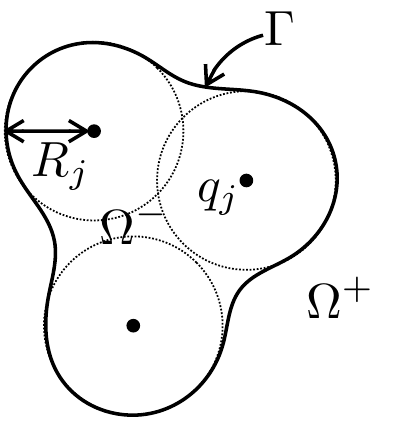}
    \caption{\label{fig:geometry} A sketch of the geometry and the interface between the molecule and solvent regions.}
\end{figure}

Electrostatics is frequently the dominant force behind molecular structure and function, and also, in the interaction with a surrounding solvent. For this reason, the development of fast and accurate models to represent electrostatics appropriately is essential to understand molecular processes. The fine-grained models in molecular dynamics describe the atomic structure explicitly and give a high level of detail, but are computationally expensive for large molecules. In contrast, continuum approximations offer an efficient alternative to compute mean-field potentials and solvation free energies \cite{RouxSimonson1999,DecherchiETal2015}. In biologically relevant settings, the solvent contains salt ions that are free to move and affect the electrostatic field. Recognizing that these ions arrange according to Boltzmann statistics at equilibrium, the electrostatic potential is well described by the Poisson-Boltzmann equation (PBE). This so-called implicit-solvent model for electrostatics considers an infinite {\it solvent} region ($\Omega^+$) with the PBE interfaced with a Poisson equation in a {\it molecule} region ($\Omega^-$) and forced by point charges at the location of the solute atoms (see Figure~\ref{fig:geometry}). Then, the resulting system of partial differential equations is given by
\begin{equation} \label{eq:poissonboltzmann}
	\begin{cases}
		- \Delta \phiext + \kappa^2 \phiext = 0, & \text{in } \Omega^+; \\
		- \Delta \phiint = \sum_{j=1}^n \frac{q_j}{\epsilonint} \delta(\mathbf{x}-\mathbf{x}_j), & \text{in } \Omega^-; \\
		\traceDe \phiext = \traceDi \phiint, & \text{at } \Gamma; \\
		\epsilonext \traceNe \phiext = \epsilonint \traceNi \phiint, & \text{at } \Gamma.
	\end{cases}
\end{equation}
Here, $\phiint$ and $\phiext$ denote the electrostatic potential in the interior molecule and exterior solvent region, respectively; $\kappa$ is the inverse of the Debye length; $n$ is the number of point sources, located at $\mathbf{x}_i$ with charge $q_i$; $\epsilonint$ and $\epsilonext$ denote the permittivity; and $\Gamma$ is the interface between $\Omega^-$ and $\Omega^+$ with unit outward pointing normal $\normal$. The operators
\begin{subequations} \label{eq:interface_conditions}
\begin{align}
    &\traceDei f(\mathbf{x}) = \lim_{\mathbf{y} \to \mathbf{x}} f(\mathbf{y}), \\
    &\traceNei f(\mathbf{x}) = \lim_{\mathbf{y} \to \mathbf{x}} \nabla f(\mathbf{y}) \cdot \normal(\mathbf{x})
\end{align}
\end{subequations}
for $\mathbf{y} \in \Omega^\pm$ and $\mathbf{x} \in \Gamma$ denote the exterior/interior Dirichlet and Neumann traces, respectively, which correspond to the limit as a point in the domain approaches the interface. In this case, the interface conditions ensure continuity of the potential ($\phi$) and electric displacement ($\epsilon\partial\phi/\partial\normal$) across $\Gamma$. This interface is enforced by the atomic radii ($R_j$ in Fig. \ref{fig:geometry}), and can be defined as the solvent-accessible\cite{LeeRichards1971}, solvent-excluded\cite{Connolly1983}, van der Waals\cite{Whitley1998}, or Gaussian surface\cite{YuETal2006}. In this work, we use the solvent-excluded surface, which avoids any geometrical singularities (i.e., cusps).

There are several numerical implementations of the PBE, spanning from finite difference \cite{BakerETal2001,JurrusETal2018,GilsonETal1985,LiETal2019} and finite element, \cite{BondETal2010,HolstETal2012} to boundary element \cite{LuETal2006,Bajaj2011,GengKrasny2013,CooperBardhanBarba2014} methods. In particular, the boundary element method (BEM) is well suited to numerically solve the PBE, especially when high precision is required \cite{CooperBardhanBarba2014}. The BEM provides an accurate description of the solute-solvent interface since it uses surface meshes. Furthermore, the Green's functions model point-charge distributions and unbounded domains precisely, and no artificial boundaries are necessary to limit the infinite solvent region. The first coupled BEM-BEM approach for the implicit-solvent model was presented by Yoon and Lenhoff \cite{YoonLenhoff1990} in 1990, to which we refer as the {\it direct} formulation, as it only couples the electrostatic potential on the molecular surface. Later, alternative boundary integral formulations were derived by Juffer {\it et al.} \cite{JufferETal1991} and Lu {\it et al.} \cite{LuETal2006}, which not only couple the potential, but also its derivative. These formulations yield a better conditioned linear system of equations and, therefore, faster convergence of the linear solver. Researchers usually choose one of these formulations, which were compared recently by Wang {\it et al.} \cite{wang2021high}, but there is no systematic study of the different formulations. Moreover, other possible BEM formulations have not been explored at all.

In this work, we took inspiration from the computational acoustics community to develop a general approach to the BEM-BEM coupling algorithm \cite{wout2021benchmarking}. Similar to the PBE, the Helmholtz equation for acoustic wave propagation can be solved efficiently with the BEM. We exploited the mathematical similarities of the acoustic Helmholtz equation and the PBE to generate a family of boundary integral formulations for the implicit-solvent model, where the Lu\cite{LuETal2006} and Juffer\cite{JufferETal1991} formulations are special cases. Also, it allowed us to design efficient formulations that are new to the PBE.

A second novelty of this study is the application of operator preconditioning to the BEM for the PBE. A large proportion of the simulation time is consumed by the iterative linear solver (e.g.~GMRES), whose convergence strongly depends on the boundary integral formulation of choice. Preconditioning can significantly reduce the solution time when it is correctly designed for the model at hand. So far, all preconditioners for the PBE are algebraic, that is, they are designed based on knowledge of the matrix. Examples are block-diagonal~\cite{altman2009accurate} and incomplete LU~\cite{chen2018preconditioning} preconditioners. Here, we use a different class of preconditioners, called operator preconditioners\cite{hiptmair2006operator,kirby2010functional}, which are designed based on knowledge of the continuous formulation. These have a strong mathematical foundation, and are straightforward to combine with the fast multipole method (FMM)~\cite{greengard1987fast}. Examples from acoustics and electromagnetics include mass\cite{betcke2020product} and Calderón\cite{andriulli2008multiplicative} preconditioning. We adapted this strategy and applied operator preconditioning to the PBE to design a novel  scaled mass-matrix preconditioning strategy. This approach incorporates knowledge on the spectrum of eigenvalues in the design of the preconditioner, and does not come with any computational overhead compared to standard preconditioners. 

The extensive analysis of formulations and preconditioners in this work allows us to give guidelines towards choosing the most efficient algorithm to solve the PBE for a given type of molecule. This is essential for our goal of designing optimal boundary integral formulations for molecular electrostatics.



The next section presents the methodology, which describes the generalized form of the BEM formulation and the preconditioning strategies. Later, the Results and Discussion section compares the computational performance of different formulations and preconditioners for a single model. Then, the performance of the most efficient preconditioned formulations are assessed on a broader set of molecules. Finally, the last section presents conclusions and an outlook for future work.

\section*{\sffamily \Large METHODOLOGY}
The implicit-solvent model for electrostatics considers a coupled system of the Poisson and Poisson-Boltzmann equations, as shown by Eq.~\eqref{eq:poissonboltzmann}. This section explains the design of boundary integral formulations to solve said system, and the development of preconditioning strategies that improve their efficiency.

\subsection*{\sffamily \large Boundary integral formulations of the Poisson-Boltzmann equation}

Since the solvent and molecule regions are assumed to be homogeneous, the Green's function of the coupled system in Eq.~\eqref{eq:poissonboltzmann} is given by
\begin{subequations}
\begin{align}
	G_\text{ext} &= \frac{e^{-\kappa ||\mathbf{x} - \mathbf{y}||}}{4\pi ||\mathbf{x} - \mathbf{y}||}, & \mathbf{x},\mathbf{y} \in \Omega^+, \mathbf{x} \ne \mathbf{y}; \\
	G_\text{int} &= \frac1{4\pi ||\mathbf{x} - \mathbf{y}||}, & \mathbf{x},\mathbf{y} \in \Omega^-, \mathbf{x} \ne \mathbf{y};
\end{align}
\end{subequations}
which correspond to the free space Green's function of the linearized Poisson-Boltzmann and Poisson equations, respectively. Applying Green's second identity to Eq.~\eqref{eq:poissonboltzmann} yields a representation of the electrostatic potential in $\Omega^+$ and $\Omega^-$ by surface potentials at the interface of the two regions. That is, 
\begin{subequations}
\label{eq:representation}
\begin{align}
    \uext(\mathbf{x}) &= \DLPext[\traceDe\uext](\mathbf{x}) - \SLPext[\traceNe\uext](\mathbf{x}), && \mathbf{x} \in \Omega^+; \label{eq:representation:ext} \\
    \uint(\mathbf{x}) &= \SLPint[\traceNi\uint](\mathbf{x}) - \DLPint[\traceDi\uint](\mathbf{x}) + \sum_{j=1}^n \frac{q_j}{\epsilonint} \frac1{4\pi||\mathbf{x}-\mathbf{x}_j||}, && \mathbf{x} \in \Omega^-; \label{eq:representation:int}
\end{align}
\end{subequations}
where
\begin{subequations}
\begin{align}
    \SLPei[\psi](\mathbf{x}) &= \iint_\Gamma G_\mathrm{ext,int}(\mathbf{x},\mathbf{y}) \psi(\mathbf{y}) \dy, && \mathbf{x} \in \Omega^\pm; \\
    \DLPei[\phi](\mathbf{x}) &= \iint_\Gamma \frac{\partial G_\mathrm{ext,int}(\mathbf{x},\mathbf{y})}{\partial\normal} \phi(\mathbf{y}) \dy, && \mathbf{x} \in \Omega^\pm;
\end{align}
\end{subequations}
are the single-layer and double-layer potential operators for each domain\cite{nedelec2001acoustic, steinbach2008numerical, sauter2010boundary}.

A common quantity of interest is the polar component of the solvation free energy. This free energy is equivalent to the work required to dissolve a molecule, and can be separated into non-polar and polar components. The non-polar part corresponds to the energy spent generating a solute-shaped cavity in the solvent, whereas the polar part is responsible for the process of charging up the cavity with the point charges $q_j$. The polar part of the solvation free energy is\cite{Baker2004}
\begin{equation}\label{eq:dG}
\Delta G = \frac{1}{2}\int_{\Omega^-} \rho(\mathbf{x})\phi_\text{reac}(\mathbf{x})\dx = \frac{1}{2} \sum_{j=1}^{n} q_j\phi_\text{reac}(\mathbf{x}_j)
\end{equation}
where 
\begin{equation}\label{eq:phireac}
\phi_\text{reac}(\mathbf{x})=\uint(\mathbf{x}) - \sum_{j=1}^n \frac{q_j}{\epsilonint} \frac1{4\pi||\mathbf{x}-\mathbf{x}_j||}
\end{equation}
is the so-called {\it reaction potential}, that arises due to the polarization of the surrounding solvent and redistribution of the salt ions. The reaction potential can be computed by subtracting out the Coulomb contribution from Eq.~\eqref{eq:representation:int}, leaving
\begin{equation}\label{eq:phireac_bie}
\phi_\text{reac}(\mathbf{x}) = \SLPint[\traceNi\uint](\mathbf{x}) - \DLPint[\traceDi\uint](\mathbf{x}).
\end{equation}
Eq.~\eqref{eq:phireac_bie} shows that the solvation free energy can be computed directly from the electrostatic potential at the surface of the molecule, which is, so far, unknown. We calculate this surface distribution using boundary integral equations.

Note that the representation formulas in Eq.~\eqref{eq:representation} relate the electrostatic potential in the exterior and interior with the electrostatic potential on the surface. To obtain a consistent formulation, a limiting process towards the surface needs to be applied to the representation formulas. Since the interface~$\Gamma$ is a smooth surface, the traces of the potential operators satisfy the so-called jump relations\cite{steinbach2008numerical}:
\begin{subequations}
\begin{align}
    \traceDei\left(\SLP[\psi]\right) &= \SL\psi; &
    \traceNei\left(\SLP[\psi]\right) &= \AD\psi \mp \frac12\psi; \\
    \traceDei\left(\DLP[\phi]\right) &= \DL\phi \pm \frac12\phi; &
    \traceNei\left(\DLP[\phi]\right) &= -\HS\phi;
\end{align}
\end{subequations}
where
\begin{subequations}\label{eq:layer_potentials}
\begin{align}
    \SLei[\psi](\mathbf{x}) &= \iint_\Gamma G_\mathrm{ext,int}(\mathbf{x},\mathbf{y}) \psi(\mathbf{y}) \dy, && \mathbf{x} \in \Gamma; \\
    \DLei[\phi](\mathbf{x}) &= \iint_\Gamma \frac{\partial G_\mathrm{ext,int}(\mathbf{x},\mathbf{y})}{\partial\normal(\mathbf{y})} \phi(\mathbf{y}) \dy, && \mathbf{x} \in \Gamma; \\
    \ADei[\psi](\mathbf{x}) &= \frac{\partial}{\partial\normal(\mathbf{x})} \iint_\Gamma G_\mathrm{ext,int}(\mathbf{x},\mathbf{y}) \psi(\mathbf{y}) \dy, && \mathbf{x} \in \Gamma; \\
    \HSei[\phi](\mathbf{x}) &= -\frac{\partial}{\partial\normal(\mathbf{x})} \iint_\Gamma \frac{\partial G_\mathrm{ext,int}(\mathbf{x},\mathbf{y})}{\partial\normal(\mathbf{y})} \phi(\mathbf{y}) \dy, && \mathbf{x} \in \Gamma
\end{align}
\end{subequations}
are the single-layer, double-layer, adjoint double-layer, and hypersingular boundary integral operators, respectively. 
Note that the evaluation points ($\mathbf{x}$) in Eq.~\eqref{eq:layer_potentials} are now on the interface $\Gamma$.
The traces of the representation formulas in Eq.~\eqref{eq:representation} are given by
\begin{subequations}
\begin{align*}
    \traceDe\uext(\mathbf{x}) &= \DLext[\traceDe\uext](\mathbf{x}) + \frac12\traceDe\uext(\mathbf{x}) - \SLext[\traceNe\uext](\mathbf{x}), \\
    \traceNe\uext(\mathbf{x}) &= -\HSext[\traceDe\uext](\mathbf{x}) - \ADext[\traceNe\uext](\mathbf{x}) + \frac12\traceNe\uext(\mathbf{x}), \\
    \traceDi\uint(\mathbf{x}) &= \SLint[\traceNi\uint](\mathbf{x}) - \DLint[\traceDi\uint](\mathbf{x}) + \frac12\traceDi\uint(\mathbf{x}) + \sum_{j=1}^n \frac{q_j}{\epsilonint} \frac1{4\pi||\mathbf{x}-\mathbf{x}_j||}, \\
    \traceNi\uint(\mathbf{x}) &= \ADint[\traceNi\uint](\mathbf{x}) + \frac12\traceNi\uint(\mathbf{x}) + \HSint[\traceDi\uint](\mathbf{x}) + \sum_{i=1}^n \frac{q_i}{\epsilonint} \frac\partial{\partial\normal(\mathbf{x})} \frac1{4\pi||\mathbf{x}-\mathbf{x}_j||},
\end{align*}
\end{subequations}
for $\mathbf{x} \in \Gamma$.
These four equations can be summarized as
\begin{subequations}
\label{eq:bie}
\begin{align}
	\begin{bmatrix} \tfrac12 \ID - \DLext & \SLext \\ \HSext & \tfrac12 \ID + \ADext \end{bmatrix}
		\begin{bmatrix} \traceDe \uext \\ \traceNe \uext \end{bmatrix}
		&= \begin{bmatrix} 0 \\ 0 \end{bmatrix}, \label{eq:bie:ext:ext} \\
	\begin{bmatrix} \tfrac12 \ID + \DLint & -\SLint \\ -\HSint & \tfrac12 \ID - \ADint \end{bmatrix}
		\begin{bmatrix} \traceDi \uint \\ \traceNi \uint \end{bmatrix}
		&= \begin{bmatrix} \varphi_D \\ \varphi_N \end{bmatrix} \label{eq:bie:int:int}
\end{align}
\end{subequations}
where $\ID$ denotes the identity operator, and
\begin{subequations}
\label{eq:sources}
\begin{align}
    \varphi_D &= \sum_{j=1}^n \frac{q_j}{\epsilonint} \frac1{4\pi||\mathbf{x}-\mathbf{x}_j||}, \\
    \varphi_N &= \sum_{j=1}^n \frac{q_j}{\epsilonint} \frac\partial{\partial\normal(\mathbf{x})} \frac1{4\pi||\mathbf{x}-\mathbf{x}_j||}
\end{align}
\end{subequations}
the projections of the sources onto the surface.
These two independent sets of boundary integral equations for the exterior and interior domains can be coupled using the interface conditions in Eq.~\eqref{eq:poissonboltzmann}. This is achieved by eliminating either the interior or the exterior unknowns. Firstly, the exterior boundary integral formulation~\eqref{eq:bie:ext:ext} can be written as 
\begin{align}
	\begin{bmatrix} \ID & 0 \\ 0 & \frac\epsilonext\epsilonint \ID \end{bmatrix} \begin{bmatrix} \tfrac12 \ID - \DLext & \SLext \\ \HSext & \tfrac12 \ID + \ADext \end{bmatrix} \begin{bmatrix} \ID & 0 \\ 0 & \frac\epsilonint\epsilonext \ID \end{bmatrix} 
		\begin{bmatrix} \traceDe \uext \\ \frac\epsilonext\epsilonint \traceNe \uext \end{bmatrix}
		&= \begin{bmatrix} 0 \\ 0 \end{bmatrix}, \nonumber \\
	\begin{bmatrix} \tfrac12 \ID - \DLext & \frac\epsilonint\epsilonext \SLext \\ \frac\epsilonext\epsilonint \HSext & \tfrac12 \ID + \ADext \end{bmatrix}
		\begin{bmatrix} \traceDi \uint \\ \traceNi \uint \end{bmatrix}
		&= \begin{bmatrix} 0 \\ 0 \end{bmatrix}. \label{eq:bie:ext:int}
\end{align}
Secondly, the interior boundary integral formulation~\eqref{eq:bie:int:int} can be written as 
\begin{align}
	\begin{bmatrix} \ID & 0 \\ 0 & \frac\epsilonint\epsilonext \ID \end{bmatrix} \begin{bmatrix} \tfrac12 \ID + \DLint & -\SLint \\ -\HSint & \tfrac12 \ID - \ADint \end{bmatrix} \begin{bmatrix} \ID & 0 \\ 0 & \frac\epsilonext\epsilonint \ID \end{bmatrix}
		\begin{bmatrix} \traceDi \uint \\ \frac\epsilonint\epsilonext \traceNi \uint \end{bmatrix}
		&= \begin{bmatrix} \ID & 0 \\ 0 & \frac\epsilonint\epsilonext \ID \end{bmatrix} \begin{bmatrix} \varphi_D \\ \varphi_N \end{bmatrix}, \nonumber \\
	\begin{bmatrix} \tfrac12 \ID + \DLint & -\frac\epsilonext\epsilonint \SLint \\ -\frac\epsilonint\epsilonext \HSint & \tfrac12 \ID - \ADint \end{bmatrix}
		\begin{bmatrix} \traceDe \uext \\ \traceNe \uext \end{bmatrix}
		&= \begin{bmatrix} \varphi_D \\ \frac\epsilonint\epsilonext \varphi_N \end{bmatrix}.\label{eq:bie:int:ext}
\end{align}
Note that four independent boundary integral equations are present while the four unknowns are reduced to two unknowns with the interface conditions. We can take linear combinations of the four boundary integral equations\cite{wout2021benchmarking} in Eq.~\eqref{eq:bie}, to produce different boundary integral formulations for the same Poisson-Boltzmann system\cite{bardhan2009numerical}.

\subsubsection*{\sffamily \normalsize The direct formulation}

The easiest choice is to take the Dirichlet traces only. Taking the first row in Eq.~\eqref{eq:bie:int:int} and the first row in Eq.~\eqref{eq:bie:ext:int} results in
\begin{align}
	\begin{bmatrix} \tfrac12 \ID + \DLint & -\SLint \\ \tfrac12 \ID - \DLext & \frac\epsilonint\epsilonext \SLext \end{bmatrix}
	\begin{bmatrix} \traceDi \uint \\ \traceNi \uint \end{bmatrix}
	&= \begin{bmatrix} \varphi_D \\ 0 \end{bmatrix}
	\label{eq:bie:dirichlet:int}
\end{align}
called the \emph{internal direct} formulation\cite{altman2009accurate}.
On the other hand, taking the first row in Eq.~\eqref{eq:bie:ext:ext} and the first row in Eq.~\eqref{eq:bie:int:ext} results in
\begin{align}
	\begin{bmatrix} \tfrac12 \ID - \DLext & \SLext \\ \tfrac12 \ID + \DLint & -\frac\epsilonext\epsilonint \SLint \end{bmatrix}
	\begin{bmatrix} \traceDe \uext \\ \traceNe \uext \end{bmatrix}
	&= \begin{bmatrix} 0 \\ \varphi_D \end{bmatrix}
	\label{eq:bie:dirichlet:ext}
\end{align}
and we obtain the \emph{external direct} formulation\cite{YoonLenhoff1990}. Notice that the order of the two equations is arbitrary and each row can be permuted to obtain an equivalent model. These are called the \emph{permuted} versions in the computational benchmarks.

The design of the direct formulation can be simplified by applying Green's identities directly on the electrostatic potential~\cite{bardhan2009numerical}. The reason to include the normal derivative of the electrostatic potential in this study is that it yields a more general design framework, where the electrostatic potential and its normal derivative at the interface can be combined to achieve better conditioned formulations.

\subsubsection*{\sffamily \normalsize The combined field formulation}

Until now, the literature on boundary integral formulations for the PBE follows specific choices to combine the electrostatic potential and its normal derivative at the surface~\cite{bardhan2009numerical}. Here, we generalize this approach by introducing the  arbitrary constants $\alpha$ and $\beta$, to take linear combinations of the boundary integral formulations in Eqs.~\eqref{eq:bie:int:int} and~\eqref{eq:bie:ext:int}. This leads to the so-called \emph{combined field integral equations} (CFIE)\cite{mitzner1966acoustic, mautz1977electromagnetic}, where
\begin{align}
	\left( \begin{bmatrix} \tfrac12 \ID + \DLint & -\SLint \\ -\HSint & \tfrac12 \ID -\ADint \end{bmatrix} + \begin{bmatrix} \alpha & 0 \\ 0 & \beta \end{bmatrix} \begin{bmatrix} \tfrac12 \ID - \DLext & \frac\epsilonint\epsilonext \SLext \\ \frac\epsilonext\epsilonint \HSext & \tfrac12 \ID + \ADext \end{bmatrix} \right)
		\begin{bmatrix} \traceDi \uint \\ \traceNi \uint \end{bmatrix}
		&= \begin{bmatrix} \varphi_D \\ \varphi_N \end{bmatrix}
		\label{eq:bie:combined:int}
\end{align}
is the internal CFIE. Alternatively, a linear combination of Eq.~\eqref{eq:bie:ext:ext} and~\eqref{eq:bie:int:ext} yields
\begin{align}
	\left( \begin{bmatrix} \tfrac12 \ID - \DLext & \SLext \\ \HSext & \tfrac12 \ID + \ADext \end{bmatrix} + \begin{bmatrix} \alpha & 0 \\ 0 & \beta \end{bmatrix} \begin{bmatrix} \tfrac12 \ID + \DLint & -\frac\epsilonext\epsilonint \SLint \\ -\frac\epsilonint\epsilonext \HSint & \tfrac12 \ID - \ADint \end{bmatrix} \right)
		\begin{bmatrix} \traceDe \uext \\ \traceNe \uext \end{bmatrix}
		&= \begin{bmatrix} \alpha & 0 \\ 0 & \beta \end{bmatrix} \begin{bmatrix} \varphi_D \\ \frac\epsilonint\epsilonext \varphi_N \end{bmatrix}
		\label{eq:bie:combined:ext}
\end{align}
which is the external CFIE.
Even though $\alpha$ and $\beta$ could be complex-valued or even boundary integral operators, we restrict ourselves to real values. Here, the \emph{internal} and \emph{external} versions of the CFIE concern the unknown potentials in the solution vector. In both cases, the entire PBE for the molecule and solvent regions are solved.

Below, we will consider specific choices of the parameters $\alpha$ and $\beta$ to obtain existing and novel formulations to solve the PBE.

{\sffamily \small The PMCHWT formulation}\\

Let us consider the choice of $\alpha = -1$ and $\beta = -1$. Then, Eq.~\eqref{eq:bie:combined:int} reads
\begin{align}
	\begin{bmatrix} -\DLint - \DLext & \SLint + \frac\epsilonint\epsilonext \SLext \\ \HSint + \frac\epsilonext\epsilonint \HSext & \ADint + \ADext \end{bmatrix}
	\begin{bmatrix} \traceDi \uint \\ \traceNi \uint \end{bmatrix}
	&= -\begin{bmatrix} \varphi_D \\ \varphi_N \end{bmatrix}
	\label{eq:bie:pmchwt:int}
\end{align}
and Eq.~\eqref{eq:bie:combined:ext} reads
\begin{align}
	\begin{bmatrix} -\DLext - \DLint & \SLext + \frac\epsilonext\epsilonint \SLint \\ \HSext + \frac\epsilonint\epsilonext \HSint & \ADext + \ADint \end{bmatrix}
	\begin{bmatrix} \traceDe \uext \\ \traceNe \uext \end{bmatrix}
	&= -\begin{bmatrix} \varphi_D \\ \frac\epsilonint\epsilonext \varphi_N \end{bmatrix}.
	\label{eq:bie:pmchwt:ext}
\end{align}
These are Fredholm equations of the first kind and known as the Poggio-Miller-Chang-Harrington-Wu-Tsai (PMCHWT) formulation in the electromagnetics and acoustics communities~\cite{poggio1973integral, chang1974surface, wu1977scattering-bor}.

{\sffamily \small The M\"uller formulation}\\

Let us consider the choice of $\alpha = 1$ and $\beta = 1$. Then, Eq.~\eqref{eq:bie:combined:int} reads
\begin{align}
	\begin{bmatrix} \ID + \DLint - \DLext & -\SLint + \frac\epsilonint\epsilonext \SLext \\ -\HSint + \frac\epsilonext\epsilonint \HSext & \ID - \ADint + \ADext \end{bmatrix}
	\begin{bmatrix} \traceDi \uint \\ \traceNi \uint \end{bmatrix}
	&= \begin{bmatrix} \varphi_D \\ \varphi_N \end{bmatrix}
	\label{eq:bie:muller:int}
\end{align}
and Eq.~\eqref{eq:bie:combined:ext} reads
\begin{align}
	\begin{bmatrix} \ID - \DLext + \DLint & \SLext - \frac\epsilonext\epsilonint \SLint \\ \HSext - \frac\epsilonint\epsilonext \HSint & \ID + \ADext - \ADint \end{bmatrix}
	\begin{bmatrix} \traceDe \uext \\ \traceNe \uext \end{bmatrix}
	&= \begin{bmatrix} \varphi_D \\ \frac\epsilonint\epsilonext \varphi_N \end{bmatrix}.
	\label{eq:bie:muller:ext}
\end{align}
These are Fredholm equations of the second kind, and their variants for electromagnetics and acoustics are known as the Müller formulation~\cite{muller1957grundprobleme}.

{\sffamily \small The Juffer formulation}\\

Let us consider the choice of $\alpha = \frac\epsilonext\epsilonint$ and $\beta = \frac\epsilonint\epsilonext$. Then, Eq.~\eqref{eq:bie:combined:int} reads
\begin{align}
	\begin{bmatrix} \tfrac12 (1 + \frac\epsilonext\epsilonint) \ID + \DLint - \frac\epsilonext\epsilonint \DLext& -\SLint + \SLext \\ -\HSint + \HSext & \tfrac12 (1 + \frac\epsilonint\epsilonext) \ID - \ADint + \frac\epsilonint\epsilonext \ADext \end{bmatrix}
		\begin{bmatrix} \traceDi \uint \\ \traceNi \uint \end{bmatrix}
		&= \begin{bmatrix} \varphi_D \\ \varphi_N \end{bmatrix}
		\label{eq:bie:juffer:int}
\end{align}
which is a Fredholm equation of the second kind, and is known as the Juffer formulation\cite{JufferETal1991}.

{\sffamily \small The Lu formulation}\\

Let us consider the choice of $\alpha = \frac\epsilonint\epsilonext$ and $\beta = \frac\epsilonext\epsilonint$. Then, Eq.~\eqref{eq:bie:combined:ext} reads
\begin{align*}
	\begin{bmatrix} \tfrac12 (1 + \frac\epsilonint\epsilonext) \ID - \DLext + \frac\epsilonint\epsilonext \DLint & \SLext - \SLint \\ \HSext - \HSint & \tfrac12 (1 + \frac\epsilonext\epsilonint) \ID + \ADext - \frac\epsilonext\epsilonint \ADint \end{bmatrix}
		\begin{bmatrix} \traceDe \uext \\ \traceNe \uext \end{bmatrix}
		&= \begin{bmatrix} \frac\epsilonint\epsilonext \varphi_D \\ \varphi_N \end{bmatrix}
\end{align*}
and multiplying the second row by $\frac\epsilonint\epsilonext$ yields
\begin{align}
	\begin{bmatrix} \tfrac12 (1 + \frac\epsilonint\epsilonext) \ID - \DLext + \frac\epsilonint\epsilonext \DLint & \SLext - \SLint \\ \frac\epsilonint\epsilonext \HSext - \frac\epsilonint\epsilonext \HSint & \tfrac12 (1 + \frac\epsilonint\epsilonext) \ID + \frac\epsilonint\epsilonext \ADext - \ADint \end{bmatrix}
		\begin{bmatrix} \traceDe \uext \\ \traceNe \uext \end{bmatrix}
		&= \begin{bmatrix} \frac\epsilonint\epsilonext \varphi_D \\ \frac\epsilonint\epsilonext \varphi_N \end{bmatrix}
		\label{eq:bie:lu:ext}
\end{align}
which is a Fredholm equation of the second kind, and we call it the Lu formulation\cite{LuETal2006, lu2009adaptive}.

\subsection*{\sffamily \large Numerical discretization}

The boundary integral formulations involve continuous operators that need to be discretized to approximate the electrostatic potentials. Here, we used standard Galerkin methods for the weak formulation with continuous piecewise linear (P1) elements on a triangular surface mesh\cite{smigaj2015solving}, implemented in the Bempp-cl library~\cite{betcke2021bempp}. Hence, the number of degrees of freedom in the model is the number of vertices in the surface mesh. Furthermore, the dense matrix arithmetic was accelerated with the fast multipole method (FMM)\cite{greengard1987fast, wang2021high} for large problems.

\subsection*{\sffamily \large Preconditioning}

The simulation of the electrostatic potential of large molecules requires a computational domain with many degrees of freedom that generate big linear systems of equations. This makes direct solvers too expensive, and iterative linear solvers are necessary. Since the BEM systems are indefinite, the generalized minimum residual (GMRES) algorithm\cite{saad1986gmres} is the most common choice of solver. The convergence of GMRES depends on the condition number and eigenvalue spectrum of the system matrix. In particular, a low condition number and clustering of eigenvalues in the complex plane are beneficial for convergence\cite{antoine2021introduction}. 

The most common strategy to improve the convergence of iterative linear solvers is preconditioning. This technique is based on solving the system $P^{-1} A \mathbf{x} = P^{-1} \mathbf{b}$ (instead of $A \mathbf{x} = \mathbf{b}$) with the preconditioner $P$. Generally speaking, when $P$ is a good approximation of $A$, the number of GMRES iterations is reduced, and when $P$ is easy to invert, the computational overhead in each GMRES iteration is small. These are two competing goals and the design of effective preconditioners is highly problem dependent. Below, we distinguish two strategies: algebraic and operator preconditioning.

\subsubsection*{\sffamily \normalsize Algebraic preconditioners}

Algebraic preconditioners are based on information regarding the discrete linear system. The most common approach is to take $P$ to be a sparse approximation of $A$ and calculate an LU-factorization, which is known as incomplete LU (ILU) preconditioning\cite{saad2003iterative}. Examples of this technique include choosing $P$ to be the diagonal of the (blocks of the) matrix\cite{altman2009accurate}, or the near interactions in a treecode-type algorithm\cite{chen2018preconditioning}.  In particular, the direct formulation in Eq.~\eqref{eq:bie:dirichlet:int} has been successfully preconditioned with the inverse of the diagonals of each block\cite{altman2009accurate,CooperBardhanBarba2014}, and we use it in the present work.

These ILU preconditioners often improve the GMRES convergence with little computational overhead, however, they have important limitations. Firstly, algebraic preconditioners need parts of the system matrix in explicit form, which becomes cumbersome when accelerators like the FMM are used\cite{carpentieri2005combining}. Secondly, parameter tweaking is necessary for optimal results\cite{sakuma2008fast}. Finally, mathematically rigorous proofs of convergence gains are difficult to obtain and the computational complexity scales unfavorably at very large problem sizes\cite{antoine2008integral}.

\subsubsection*{\sffamily \normalsize Operator preconditioners}

Operator preconditioners are based on information regarding the continuous boundary integral formulation. They take into account the functional analysis of the boundary integral operators and make sure that the preconditioner maps toward the correct function spaces to guarantee well-conditioned formulations\cite{hiptmair2006operator, kirby2010functional}. Furthermore, they use knowledge of integral operator algebra (such as projection properties) to design preconditioners that yield second-kind boundary integral equations that have eigenvalue clustering. The advantages of this approach include the availability of a rigorous mathematical justification and the fact that preconditioners can be assembled separately from the system matrix. Operator preconditioning is straightforward to implement in existing software for the BEM\cite{betcke2020product}. Examples of operator preconditioners include opposite-order preconditioning\cite{steinbach1998construction}, Calderón preconditioning\cite{andriulli2008multiplicative}, and OSRC preconditioning\cite{antoine2005alternative}. Their effectiveness to coupled systems has been shown for electromagnetics\cite{cools2011calderon} and acoustics\cite{wout2021benchmarking}. However, no operator preconditioning for the PBE is known so far. Next, we explain two operator preconditioners, namely, mass-matrix and Calderón preconditioning, in more detail. 

{\sffamily \small Mass-matrix preconditioning}\\

Functional analysis of boundary integral formulations suggest that Fredholm integral equations of second kind are better conditioned than first-kind equations~\cite{sauter2010boundary}. Second-kind formulations have boundary integral operators in the form of
\begin{equation} \label{eq:prec:2kind}
	\begin{bmatrix} a\ID & 0 \\ 0 & b\ID \end{bmatrix} + C
\end{equation}
with $C$ a compact operator and $a$ and $b$ nonzero constants. The Müller, Juffer and Lu formulations (see Eqs.~\eqref{eq:bie:muller:int}--\eqref{eq:bie:lu:ext}) are examples of second-kind equations since the Calderón boundary integral operator, namely
\begin{equation} \label{eq:calderon}
	\begin{bmatrix} -\DL & \SL \\ \HS & \AD \end{bmatrix},
\end{equation}
is compact on smooth surfaces\cite{colton2013integral}, as is the case for the solvent-excluded surface. The eigenvalues of these systems cluster around two accumulation points in the complex plane, given by $a$ and $b$, which yields fast convergence of GMRES. Hence, second-kind formulations are already sufficiently well conditioned and simple preconditioners suffice. In this case we use mass-matrix preconditioning\cite{betcke2020product}.

The standard mass-matrix preconditioning uses the operator
\begin{equation}\label{eq:prec:mass}
	\begin{bmatrix} M & 0 \\ 0 & M \end{bmatrix}^{-1} \left( \begin{bmatrix} a\ID & 0 \\ 0 & b\ID \end{bmatrix} + C \right).
\end{equation}
The matrix $M$ is the so-called mass-matrix, where each element corresponds to the inner product between the basis functions of the range and dual-to-range spaces. As this is nonzero only when these functions overlap, $M$ is highly sparse, and a sparse LU factorization yields an efficient approximation of the inverse. In particular, $M$ is a diagonal matrix for piecewise-constant basis functions and sparse but not diagonal for piecewise-linear basis functions. From a continuous point of view, the matrix $M$ is the weak formulation of an identity operator that maps from the range to the dual-to-range space of the boundary integral formulation. This simple preconditioning strategy often improves the GMRES convergence with little computational overhead.

An important deficiency of mass-matrix preconditioning is that it is agnostic to the specific boundary integral operator (apart from its function spaces) and the convergence gain is often limited. In particular, it does not change the location of the eigenvalue accumulation points of the boundary integral formulation. Here, we propose to include the knowledge on accumulation points in the design of mass-matrix preconditioning and significantly improve convergence without increasing computation time. We propose a scaled mass-matrix preconditioning strategy of
\begin{equation}
	\begin{bmatrix} aM & 0 \\ 0 & bM \end{bmatrix}^{-1} \left( \begin{bmatrix} a\ID & 0 \\ 0 & b\ID \end{bmatrix} + C \right).
	\label{eq:prec:mass:scaled}
\end{equation}
It is straightforward to show that this preconditioning strategy yields an identity plus compact operator. Hence, the two accumulation points $a$ and $b$ in the original formulations are both mapped to the point~$1$ in the complex plane. Having a single accumulation point instead of two is beneficial to the convergence of the GMRES solver because the eigenvalues are more tightly clustered. Furthermore, the scaling does not incur any computational overhead compared to standard mass-matrix preconditioning.

The M\"uller formulations in Eqs.~\eqref{eq:bie:muller:int} and~\eqref{eq:bie:muller:ext}, and the Lu formulation in Eq.~\eqref{eq:bie:lu:ext}, have $a=b$, so no scaling of the mass-matrix preconditioner is necessary. In the case of the Juffer formulation in Eq.~\eqref{eq:bie:juffer:int}, we scale the mass-matrix with $a=1+\epsilonext/\epsilonint$ and $b=1+\epsilonint/\epsilonext$.

{\sffamily \small Calderón preconditioning}\\

First-kind formulations do not have any eigenvalue accumulation and are normally worse conditioned than second-kind formulations, but its solution is often more accurate\cite{yla2008analysis}. In this case, more powerful preconditioning strategies are needed. One of the best preconditioning strategies is called Calderón preconditioning, which is based on the Calderón projection property\cite{steinbach2008numerical, sauter2010boundary} stating
\begin{equation} \label{eq:prec:calderon}
	\begin{bmatrix} -\DL & \SL \\ \HS & \AD \end{bmatrix}^2 = \frac14 \begin{bmatrix} \ID & 0 \\ 0 & \ID \end{bmatrix}.
\end{equation}
Since the identity operator is trivial to invert, the property states that the square of the block Calderón operator is trivially solved as well. This suggests that preconditioning this system by itself becomes highly effective. Even though none of the boundary integral formulations include this Calderón matrix exactly, we use this idea of self-preconditioning for the PMCHWT formulation in Eqns.~\eqref{eq:bie:pmchwt:int} and~\eqref{eq:bie:pmchwt:ext}. We distinguish three different ways to apply Calderón preconditioning to the PMCHWT formulation: using the entire matrix as the preconditioner, or, using either the interior or exterior part. For the internal PMCHWT formulation~\eqref{eq:bie:pmchwt:int} this is
\begin{subequations}
\label{eq:prec:calderon:int}
\begin{align}
    \label{eq:prec:calderon:int:full}
	&\begin{bmatrix} -\DLint - \DLext & \SLint + \frac\epsilonint\epsilonext \SLext \\ \HSint + \frac\epsilonext\epsilonint \HSext & \ADint + \ADext \end{bmatrix}^2
		\begin{bmatrix} \traceDi \uint \\ \traceNi \uint \end{bmatrix}
		= -\begin{bmatrix} -\DLint - \DLext & \SLint + \frac\epsilonint\epsilonext \SLext \\ \HSint + \frac\epsilonext\epsilonint \HSext & \ADint + \ADext \end{bmatrix} \begin{bmatrix} \varphi_D \\ \varphi_N \end{bmatrix}, \\
    \label{eq:prec:calderon:int:int}
	&\begin{bmatrix} -\DLint & \SLint \\ \HSint & \ADint \end{bmatrix}
		\begin{bmatrix} -\DLint - \DLext & \SLint + \frac\epsilonint\epsilonext \SLext \\ \HSint + \frac\epsilonext\epsilonint \HSext & \ADint + \ADext \end{bmatrix}
		\begin{bmatrix} \traceDi \uint \\ \traceNi \uint \end{bmatrix}
		= -\begin{bmatrix} -\DLint & \SLint \\ \HSint & \ADint \end{bmatrix} \begin{bmatrix} \varphi_D \\ \varphi_N \end{bmatrix}, \\
    \label{eq:prec:calderon:int:ext}
	&\begin{bmatrix} -\DLext & \frac\epsilonint\epsilonext \SLext \\ \frac\epsilonext\epsilonint \HSext & \ADext \end{bmatrix}
		\begin{bmatrix} -\DLint - \DLext & \SLint + \frac\epsilonint\epsilonext \SLext \\ \HSint + \frac\epsilonext\epsilonint \HSext & \ADint + \ADext \end{bmatrix}
		\begin{bmatrix} \traceDi \uint \\ \traceNi \uint \end{bmatrix}
		= -\begin{bmatrix} -\DLext & \frac\epsilonint\epsilonext \SLext \\ \frac\epsilonext\epsilonint \HSext & \ADext \end{bmatrix} \begin{bmatrix} \varphi_D \\ \varphi_N \end{bmatrix}
\end{align}
\end{subequations}
and for the external PMCHWT formulation~\eqref{eq:bie:pmchwt:ext} this is given by
\begin{subequations}
\label{eq:prec:calderon:ext}
\begin{align}
    \label{eq:prec:calderon:ext:full}
	&\begin{bmatrix} -\DLext - \DLint & \SLext + \frac\epsilonext\epsilonint \SLint \\ \HSext + \frac\epsilonint\epsilonext \HSint & \ADext + \ADint \end{bmatrix}^2
	    \begin{bmatrix} \traceDe \uext \\ \traceNe \uext \end{bmatrix}
	    = -\begin{bmatrix} -\DLext - \DLint & \SLext + \frac\epsilonext\epsilonint \SLint \\ \HSext + \frac\epsilonint\epsilonext \HSint & \ADext + \ADint \end{bmatrix} \begin{bmatrix} \varphi_D \\ \frac\epsilonint\epsilonext \varphi_N \end{bmatrix}, \\
    \label{eq:prec:calderon:ext:ext}
	&\begin{bmatrix} -\DLext & \SLext \\ \HSext & \ADext \end{bmatrix}
		\begin{bmatrix} -\DLext - \DLint & \SLext + \frac\epsilonext\epsilonint \SLint \\ \HSext + \frac\epsilonint\epsilonext \HSint & \ADext + \ADint \end{bmatrix}
		\begin{bmatrix} \traceDe \uext \\ \traceNe \uext \end{bmatrix}
		= -\begin{bmatrix} -\DLext & \SLext \\ \HSext & \ADext \end{bmatrix} \begin{bmatrix} \varphi_D \\ \frac\epsilonint\epsilonext \varphi_N \end{bmatrix}, \\
    \label{eq:prec:calderon:ext:int}
	&\begin{bmatrix} -\DLint & \frac\epsilonext\epsilonint \SLint \\ \frac\epsilonint\epsilonext \HSint & \ADint \end{bmatrix}
		\begin{bmatrix} -\DLext - \DLint & \SLext + \frac\epsilonext\epsilonint \SLint \\ \HSext + \frac\epsilonint\epsilonext \HSint & \ADext + \ADint \end{bmatrix}
		\begin{bmatrix} \traceDe \uext \\ \traceNe \uext \end{bmatrix}
		= -\begin{bmatrix} -\DLint & \frac\epsilonext\epsilonint \SLint \\ \frac\epsilonint\epsilonext \HSint & \ADint \end{bmatrix} \begin{bmatrix} \varphi_D \\ \frac\epsilonint\epsilonext \varphi_N \end{bmatrix}.
\end{align}
\end{subequations}
The Calderón projection property~\eqref{eq:prec:calderon} does not hold exactly due to the linear combinations of interior and exterior boundary integral operators. However, a functional analysis of the operator products shows that the preconditioned formulation is a second-kind equation~\cite{antoine2008integral, yan2010comparative, cools2011calderon, wout2021highcontrast}, of the form of Eq.~\eqref{eq:prec:2kind}. Specifically, the accumulation points of the eigenvalues have been calculated explicitly for the acoustic equivalent~\cite{wout2021highcontrast} and carry over to the PBE. These points are
\begin{equation}\label{eq:precond:calderon:eigenval}
	\tilde\lambda_\text{full} = \frac12 + \frac\epsilonext{4\epsilonint} + \frac\epsilonint{4\epsilonext}
\end{equation}
for the full Calderón preconditioning, and the two points
\begin{equation}\label{eq:precond:calderon:eigenval_semi}
	\tilde\lambda_\text{half,1} = \frac14 + \frac\epsilonext{4\epsilonint}
	\quad \text{and} \quad
	\tilde\lambda_\text{half,2} = \frac14 + \frac\epsilonint{4\epsilonext}
\end{equation}
for the half Calderón preconditioning with either the interior or exterior operator. Notice that the two accumulation points for the half Calderón preconditioner are located far apart since the permittivity ratio ($\epsilonint/\epsilonext$) is large: over twenty for typical settings. Hence, the full Calderón preconditioner is expected to be much more effective than the half Calderón preconditioner. This issue can be mitigated with the same innovation as for the Juffer formulation, that is, a scaled mass-matrix preconditioner on top of the half Calderón preconditioner, looking to achieve a single accumulation point. In this case, we use $a=\tilde\lambda_\text{half,1}$ and $b=\tilde\lambda_\text{half,2}$ or vice versa in Eq.~\eqref{eq:prec:mass:scaled}. We expect Calderón preconditioning of first-kind formulations to improve the GMRES convergence substantially, but each iteration is much more expensive as the preconditioner is a dense matrix.




\section*{\sffamily \Large RESULTS AND DISCUSSION}

This section presents computational benchmarks to assess the performance of the different BEM formulations of the implicit-solvent model. The numerical simulations use a value of $\kappa=0.125$ \AA$^{-1}$ for the inverse Debye length (corresponding to 150 mM of salt in the solvent), and $\epsilonint=4$ and $\epsilonext=80$ for the interior and exterior permittivity, respectively. The algorithm was implemented in the open-source software Bempp-cl\cite{betcke2021bempp} for the BEM, which interfaces with ExaFMM\cite{wang2021exafmm,wang2021high} for the FMM, and SciPy's\cite{virtanen2020scipy} implementation of the GMRES linear solver. The tolerance of the solver was set to $10^{-5}$ without restart. Meshes were generated with the software NanoShaper~\cite{Decherchi2013a}.

All computations were performed on a workstation with two Intel\textregistered\ Xeon\textregistered\ E5-2680 v3 @ 2.50GHz CPUs with 12 physical cores each (48 threads in total) and 94~GB of RAM. Timings reported in this section are an average over 6 independent runs, and neglecting one outlier. 

\subsection*{\sffamily \large Performance of different formulations and preconditioners}

As an initial test to analyze the performance of the different formulations and preconditioners, we computed the solvation energy of a single arginine, which contains 27 atoms. This molecule is available in the Protein Data Bank\footnote{\url{https://www.rcsb.org/ligand/ARG}}. The solvent-excluded surface was discretized with 10 vertices/\AA$^2$, yielding a total of 4052 panels and 2028 vertices. We used the default Bempp-cl parameters, which sets a quadrature integration scheme of order 4. For such small tests, the FMM was not necessary, and matrices were assembled in dense form. The combinations of formulations and preconditioners used in this benchmark are detailed in Table~\ref{tab:formulations}.

\begin{table}
    \centering
    \def\arraystretch{0.9}
    \small
    \begin{tabular}{c | c c | c}
        Case N$^\circ$ & Formulation & Preconditioner & Fig.  \\
        \hline
         & \multicolumn{2}{c|}{Direct formulations} & \\
        1 & Internal (Eq. \eqref{eq:bie:dirichlet:int}) & None & \multirow{8}{*}{\ref{fig:arginine:direct}}\\
        2 & Internal (Eq. \eqref{eq:bie:dirichlet:int}) & Block diagonal\cite{altman2009accurate} & \\
        3 & External (Eq. \eqref{eq:bie:dirichlet:ext}) & None & \\
        4 & External (Eq. \eqref{eq:bie:dirichlet:ext}) & Block diagonal\cite{altman2009accurate} & \\
        & \multicolumn{2}{c|}{Direct permuted formulations} &\\
        5 & Internal (Eq. \eqref{eq:bie:dirichlet:int}) & None & \\
        6 & Internal (Eq. \eqref{eq:bie:dirichlet:int}) & Block diagonal\cite{altman2009accurate} & \\
        7 & External (Eq. \eqref{eq:bie:dirichlet:ext}) & None &  \\
        8 & External (Eq. \eqref{eq:bie:dirichlet:ext}) & Block diagonal\cite{altman2009accurate} & \\
        & & & \\
        &\multicolumn{2}{c|}{Second kind formulations} &\\
        9 & Juffer (Eq. \eqref{eq:bie:juffer:int}) & Mass matrix (Eq. \eqref{eq:prec:mass}) & \multirow{8}{*}{\ref{fig:arginine:secondkind}}\\
        10 & Juffer (Eq. \eqref{eq:bie:juffer:int}) & Scaled mass matrix (Eq. \eqref{eq:prec:mass:scaled}) & \\
        11 & Lu (Eq. \eqref{eq:bie:lu:ext}) & None & \\
        12 & Lu (Eq. \eqref{eq:bie:lu:ext}) & Mass matrix (Eq. \eqref{eq:prec:mass}) & \\
        13 & M\"uller Int. (Eq. \eqref{eq:bie:muller:int}) & Mass matrix (Eq. \eqref{eq:prec:mass}) & \\
        14 & M\"uller Ext. (Eq. \eqref{eq:bie:muller:ext}) & Mass matrix (Eq. \eqref{eq:prec:mass}) & \\
        &&&\\
        &\multicolumn{2}{c|}{PMCHWT formulations}&\\
        15 & Internal (Eq. \eqref{eq:bie:pmchwt:int}) & Mass matrix (Eq. \eqref{eq:prec:mass}) & \multirow{12}{*}{\ref{fig:arginine:calderon}}\\
        16 & Internal (Eq. \eqref{eq:bie:pmchwt:int}) & Calder\'on (Eq. \eqref{eq:prec:calderon:int:full})  & \\
        17 & Internal (Eq. \eqref{eq:bie:pmchwt:int}) & Internal Calder\'on (Eq. \eqref{eq:prec:calderon:int:int}) & \\
        18 & Internal (Eq. \eqref{eq:bie:pmchwt:int}) & External Calder\'on (Eq. \eqref{eq:prec:calderon:int:ext}) & \\
        19 & External (Eq. \eqref{eq:bie:pmchwt:ext}) & Mass matrix (Eq. \eqref{eq:prec:mass}) & \\
        20 & External (Eq. \eqref{eq:bie:pmchwt:ext}) &  Calder\'on (Eq. \eqref{eq:prec:calderon:ext:full})  & \\
        21 & External (Eq. \eqref{eq:bie:pmchwt:ext}) &  External Calder\'on (Eq. \eqref{eq:prec:calderon:ext:ext})  & \\
        22 & External (Eq. \eqref{eq:bie:pmchwt:ext}) &  Internal Calder\'on (Eq. \eqref{eq:prec:calderon:ext:int}) & \\
        23 & Internal (Eq. \eqref{eq:bie:pmchwt:int}) & Scaled-mass Interior Calder\'on (Eq. \eqref{eq:prec:calderon:int:int}) & \\
        24 & Internal (Eq. \eqref{eq:bie:pmchwt:int}) & Scaled-mass Exterior Calder\'on (Eq. \eqref{eq:prec:calderon:int:ext}) & \\
        25 & External (Eq. \eqref{eq:bie:pmchwt:ext}) &  Scaled-mass Exterior Calder\'on (Eq. \eqref{eq:prec:calderon:ext:ext})  & \\
        26 & External (Eq. \eqref{eq:bie:pmchwt:ext}) &  Scaled-mass Interior Calder\'on (Eq. \eqref{eq:prec:calderon:ext:int}) & \\
        \hline
    \end{tabular}
    \caption{The formulations and preconditioners used for the arginine benchmark.}
    \label{tab:formulations}
\end{table}

\begin{figure}[!ht]
    \centering
    \includegraphics[width=\columnwidth,keepaspectratio=true]{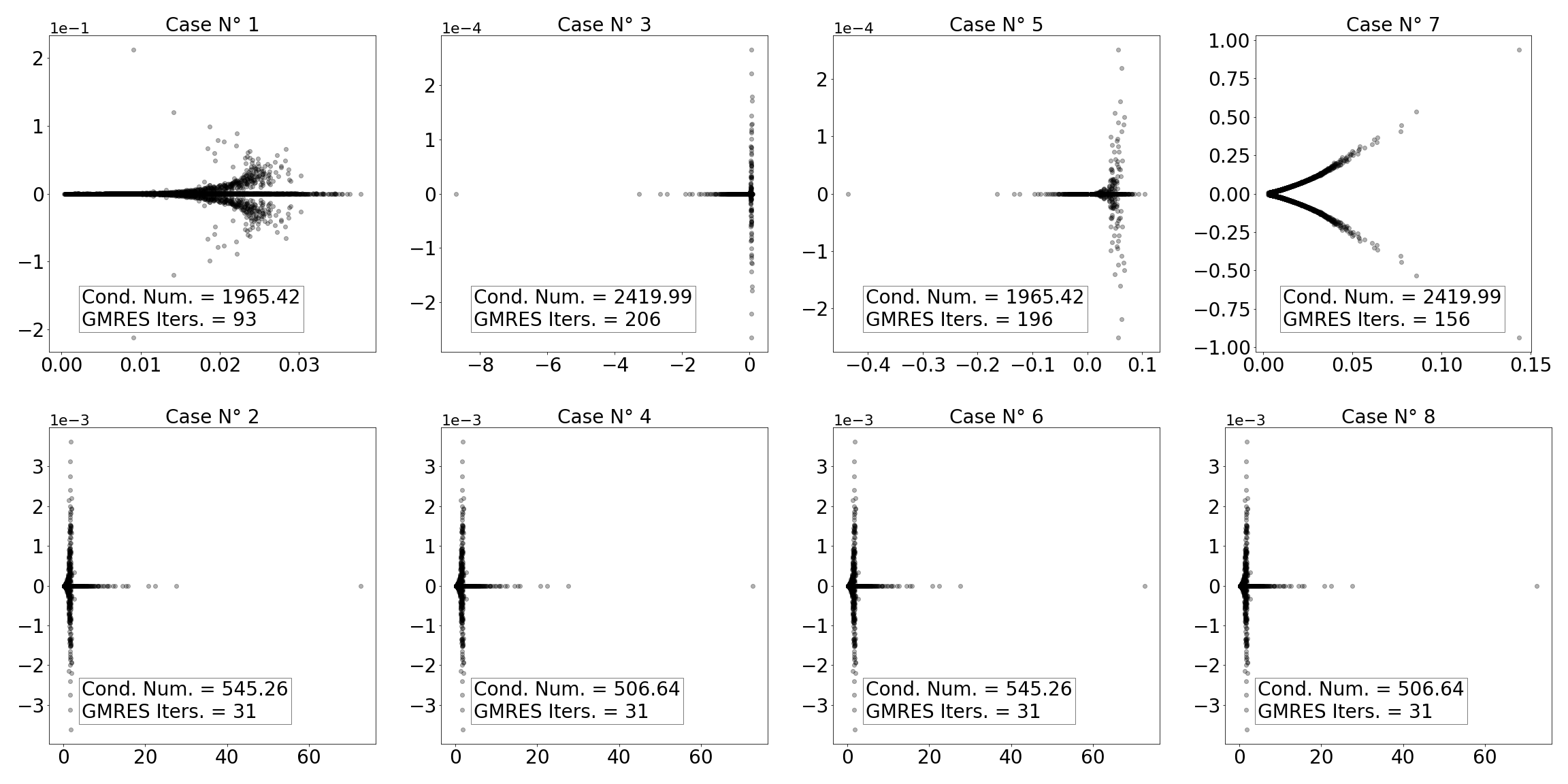}
    \caption{The spectrum of the system matrix of the direct formulations for arginine, with (bottom row) and without (top row) the block diagonal preconditioner. Each marker represents the real and imaginary part of one of the 2028 eigenvalues of the (preconditioned) system.}
    \label{fig:arginine:direct}
\end{figure}

\begin{figure}[!ht]
    \centering
    \includegraphics[width=\columnwidth,keepaspectratio=true]{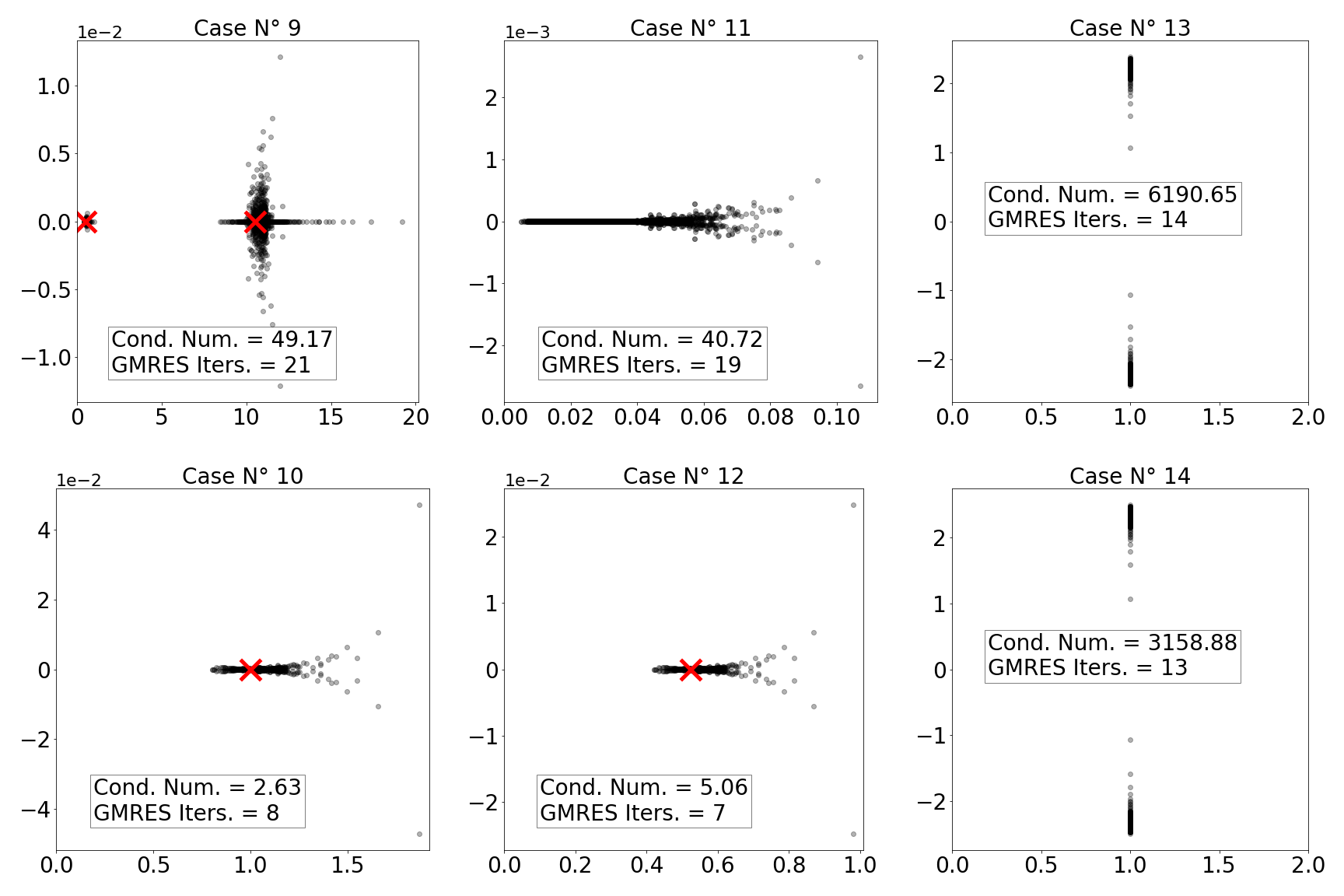}
    \caption{The spectrum of the system matrix of the second-kind formulations and (scaled) mass preconditioning for arginine. Each marker represents the real and imaginary part of one of the 2028 eigenvalues of the (preconditioned) system. The red crosses are the predicted accumulation points.}
    \label{fig:arginine:secondkind}
\end{figure}

\begin{figure}[!ht]
    \centering
    \includegraphics[width=\columnwidth,keepaspectratio=true]{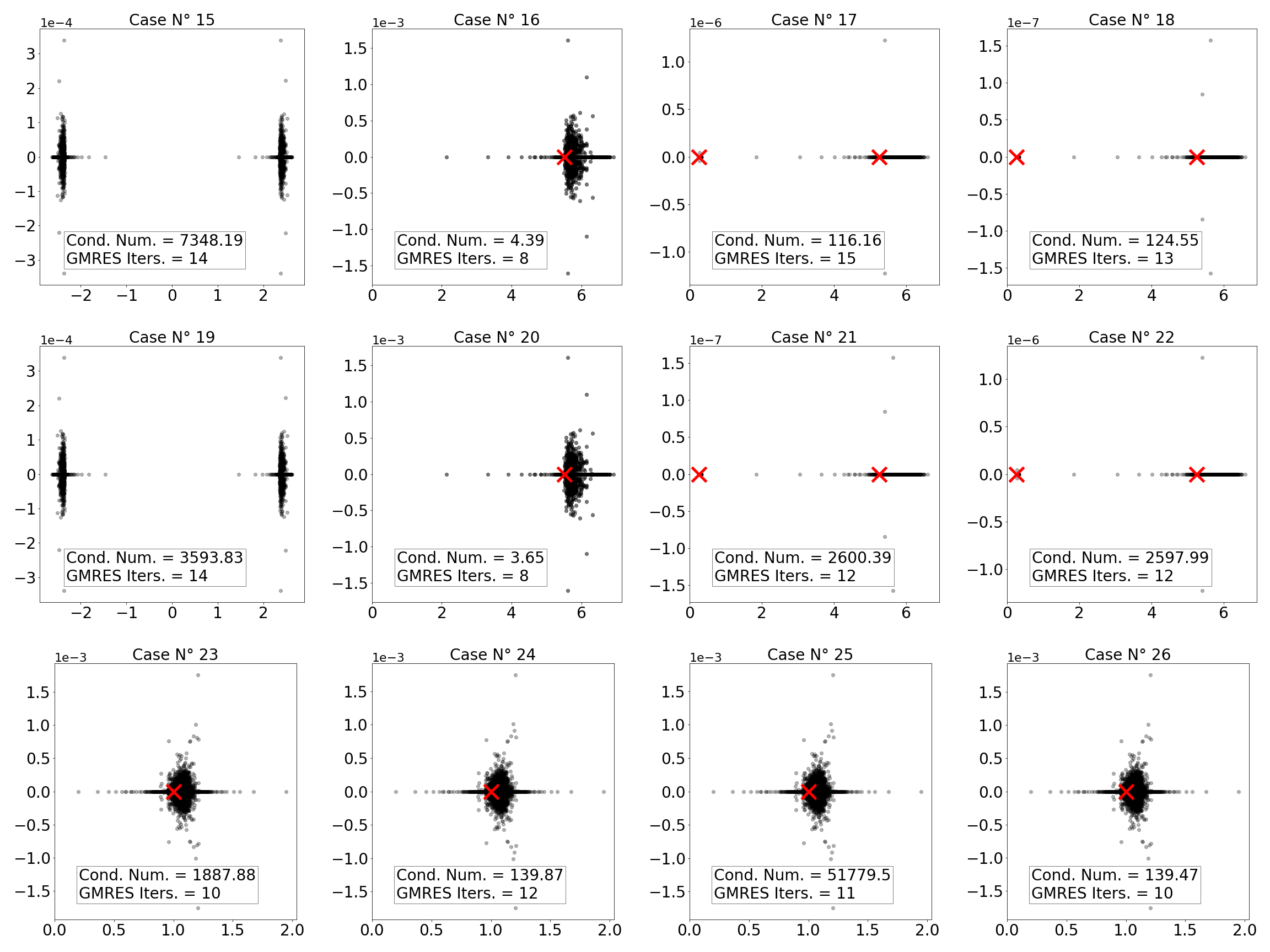}
    \caption{The spectrum of the system matrix of the first-kind formulation and Calderón preconditioning for arginine. Each marker represents the real and imaginary part of one of the 2028 eigenvalues of the (preconditioned) system. The red crosses are the predicted accumulation points.}
    \label{fig:arginine:calderon}
\end{figure}

The convergence of GMRES is often the defining parameter concerning overall efficiency of the implicit-solvent model. The number of iterations strongly depends on the eigenvalues of the linear system. As a rough measure, a small condition number yields a sharper bound on the number of iterations, however, a more precise indication of actual convergence is the clustering of the eigenvalues in the complex plane\cite{antoine2021introduction}. The goal of preconditioning is to improve the spectral properties of the system matrix without increasing the computational overhead too much. Figures \ref{fig:arginine:direct}, \ref{fig:arginine:secondkind}, and \ref{fig:arginine:calderon} show the eigenvalue spectrum, condition number and number of GMRES iterations of the arginine molecule for all models presented in Table~\ref{tab:formulations}. The calculated solvation energies for all of these formulations lay within the range between -64.137 kcal/mol and -64.278 kcal/mol.

Figure~\ref{fig:arginine:direct} confirms the improved performance of the GMRES with block-diagonal preconditioning\cite{altman2009accurate} for the direct formulation. The condition number reduces by an order of magnitude, and iteration count decreases up to 6.6~times. Comparing the spectrum from the top (no preconditioner) and bottom (block-diagonal preconditioner) rows in Figure \ref{fig:arginine:direct}, the block-diagonal preconditioner effectively clusters most of the eigenvalues, however, a significant number of outliers are present that deteriorate the conditioning of the system. The still relatively high condition numbers are thus not surprising. The direct formulations with block-diagonal preconditioner are low cost models compared to the combined formulations with operator preconditioning, still making it an attractive alternative. More specifically, the direct formulation uses only half of the boundary operators, and the block-diagonal preconditioner can be computed analytically. 

The combined boundary integral formulations involve all boundary integral operators of the Calderón operator (single, double, adjoint, and hypersingular operators), making them computationally more expensive than the direct formulation. On the upside, they are better conditioned, converge quicker, and are more stable with respect to mesh sizes. These characteristics are confirmed by the computational benchmarks in Figures~\ref{fig:arginine:secondkind} and~\ref{fig:arginine:calderon}. The combined formulations require less GMRES iterations than the direct formulations and the spectrum has clustering of eigenvalues around accumulation points. Specifically, the Juffer formulation with standard mass preconditioning (Case~9) has two accumulation points at $\frac12(1+\epsilonext/\epsilonint) = 10.5$ and $\frac12(1+\epsilonint/\epsilonext) = 0.525$, which correspond to the coefficients multiplying the identity in Eq.~\eqref{eq:bie:juffer:int}. This serves as motivation to use a scaled mass preconditioner (Case~10) as shown in Eq.~\eqref{eq:prec:mass:scaled}, which effectively joins these two accumulation points at the point~1 in the complex plane. This clustering of eigenvalues reduces the condition number and improves GMRES convergence. Alternatively, the Lu formulation already has a single accumulation point at $\frac12(1+\epsilonint/\epsilonext) = 0.525$ with standard mass preconditioning (Case~12), and has similar spectral properties. This is expected, as the coefficient next to the diagonal is the same in the top-left and bottom-right blocks of Eq.~\eqref{eq:bie:lu:ext}. Comparing Cases~11 and 12 is useful to see the impact of the simple mass-matrix preconditioner: the eigenvalues are scaled up 20$\times$, moving away from 0, and the cluster is more compact, decreasing the number of iterations by almost 3$\times$. In the M\"uller formulation of Eqs.~\eqref{eq:bie:muller:int} and \eqref{eq:bie:muller:ext} the identity is not scaled, and results on the eigenvalues largely clustered around 1 in Cases~13 and 14, however, these are scattered in the imaginary axis, leading to a high condition number.

The first-kind formulation does not have any eigenvalue clustering and is expected to be ill-conditioned, which can be solved by employing Calderón preconditioning. The results in Figure~\ref{fig:arginine:calderon} confirm this expected behavior: the eigenvalues of the mass-matrix preconditioned system (Cases~15 and 19) are clustered around two points (one positive and one negative), while the full Calderón preconditioned matrix (Cases~16 and 20) have all eigenvalues in the positive half-plane with accumulation at the expected point of $\frac12 + \frac\epsilonext{4\epsilonint} + \frac\epsilonint{4\epsilonext} = 5.5125$ (see Eq.~\eqref{eq:precond:calderon:eigenval}). Furthermore, the number of GMRES iterations are reduced by a factor of almost two. The idea of using either the interior or exterior Calderón operator as preconditioner (Cases~17, 18, 21 and 22) reduces the computational overhead but the two accumulation points $\frac14 + \frac\epsilonext{4\epsilonint} = 5.25$ and $\frac14 + \frac\epsilonint{4\epsilonext} = 0.2625$ (see Eq.~\eqref{eq:precond:calderon:eigenval_semi}), are far separated due to the large contrast in permittivity between the molecule and solvent regions. Hence, no effective eigenvalue accumulation is achieved to improve matrix conditioning. We can overcome this situation by using a scaled mass-matrix preconditioner in Cases~23 through 26. The last row in Figure~\ref{fig:arginine:calderon} shows that this technique effectively reduces accumulation to a single point. However, the GMRES iteration count decreases only slightly due to the outliers in the spectrum.

These benchmarks confirm the theory that eigenvalue clustering improves GMRES convergence~\cite{antoine2021introduction} and that operator preconditioning places the eigenvalues close to the expected accumulation point, as pointed out in the Methodology section. Since the number of GMRES iterations are significantly reduced, the benchmark also confirms that our preconditioning strategy is effective: we first analyze the accumulation points of the eigenvalues, and then design scaled mass-matrix preconditioners that accumulate all eigenvalues of the preconditioned matrix around a single point in the complex plane. This strategy works for any second-kind boundary integral formulation, including the Juffer and Lu formulations. Alternatively, first-kind formulations can be improved through Calderón preconditioning. Given that all operator preconditioned formulations have a single accumulation point, further optimizations of the conditioning should influence the spreading of the eigenvalues. Unfortunately, to the best of our knowledge, no mature theory is available to this end, and such study would need to be guided by computational benchmarks.

\subsection*{\sffamily \large Scaling with protein size}

\begin{figure}[!ht]
    \centering
    \includegraphics[width=0.95\columnwidth,keepaspectratio=true]{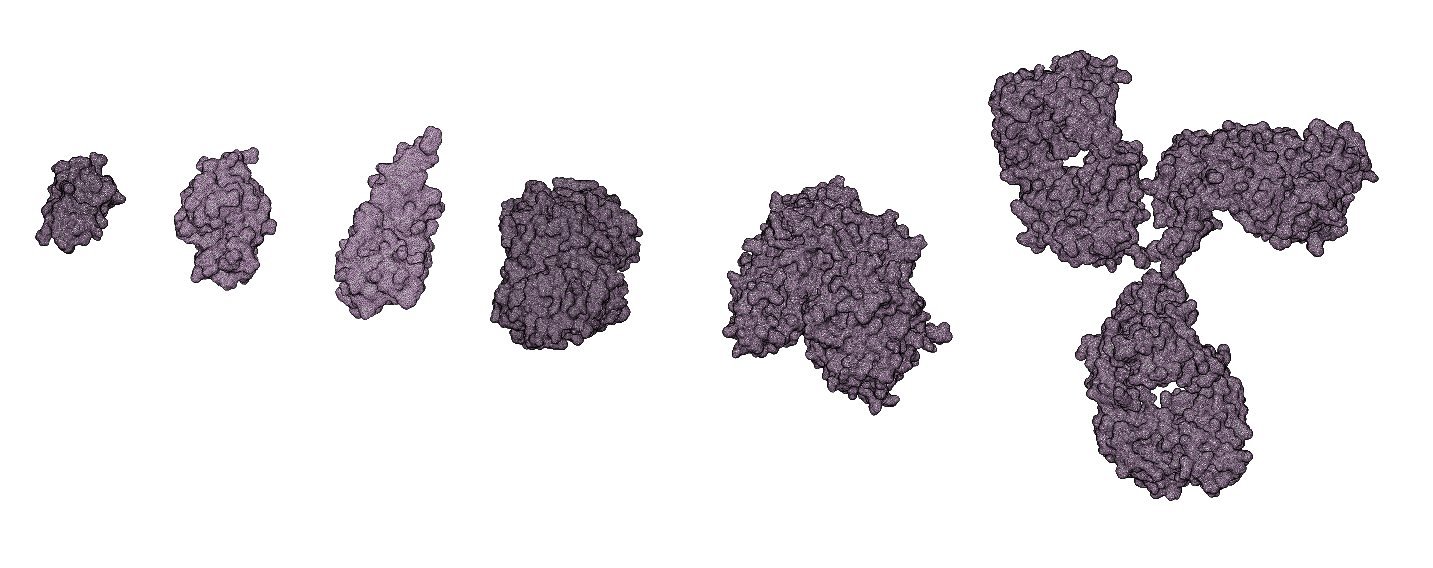}
    \caption{\label{fig:molecules} The computational domains that model the molecules (from left to right) 1bpi, 1lyz, 1a7m, 1x1z, 4lgp, and 1igt. The size of each grid is according to scale.}
\end{figure}

The benchmarks presented above confirmed the significant impact of the  boundary integral formulation and preconditioner on the computational performance of the implicit-solvent model. This choice is even more critical as the problem size increases. The computational footprint of the BEM with dense matrices scales at least quadratically with the size of the mesh, however, we can reach linear scaling through the FMM. The use of the FMM may change the optimal choice of formulation and preconditioner compared to the arginine benchmarks, where matrices were dense. Here, we study the complexity of the model with the FMM,  considering five of the best performing Cases from Table \ref{tab:formulations}, on increasingly large molecules. Specifically, the following six molecules are considered: 1bpi (12,524 vertices), 1lyz (21,844 vertices), 1a7m (31,042 vertices), 1x1z (55,446 vertices), 4lgp (116,062 vertices), and 1igt (226,312 vertices), which are shown in in Figure~\ref{fig:molecules}. The number of vertices correspond to the size of the triangular surface mesh with a density of 4 vertices/\AA$^2$. We parameterized these molecules with the PARSE force field using PDB2PQR~\cite{Dolinsky04}. As preconditioned formulations, let us consider:
\begin{itemize}
    \item The internal direct formulation with block diagonal preconditioning (Case 2 in Table~\ref{tab:formulations});
    \item The Juffer formulation with scaled mass matrix preconditioning (Case 10 in Table~\ref{tab:formulations});
    \item The Lu formulation with mass matrix preconditioning (Case 12 in Table~\ref{tab:formulations});
    \item The external PMCHWT formulation with full Calderón preconditioning (Case 20 in Table~\ref{tab:formulations}).
    \item The internal PMCHWT formulation with scaled-mass interior Calderón preconditioning (Case 23 in Table~\ref{tab:formulations}).
\end{itemize}
These models were chosen because they are the best performing versions within their respective families of boundary integral formulations.

Table~\ref{tab:results} contains the results of these tests. The FMM was used for all benchmarks, with a maximum of 50 Gauss points per leaf box (`ncrit') and an expansion order of 3. The Gauss quadrature order was set to 4 for the singular and 3 for the regular parts of the near interactions. The FMM was also applied to the calculation of the right-hand-side (RHS) vector, with an `ncrit' of 500 and an expansion order of 10. We used a high expansion order  to ensure good accuracy in the RHS, moreover, this calculation time was still small compared to the rest of the algorithm.

\begin{table}
\centering
\resizebox{\textwidth}{!}{%
\def\arraystretch{0.8}
\begin{tabular}{crrrrrrrrrr}
\hline
\multirow{2}{*}{\def\arraystretch{0.7}\begin{tabular}[c]{@{}c@{}}Case\\ N°\end{tabular}} & \multicolumn{6}{c}{Timings (s)} & \multicolumn{1}{c}{\multirow{2}{*}{\def\arraystretch{0.7}\begin{tabular}[c]{@{}c@{}}GMRES\\ iterations\end{tabular}}} & \multicolumn{1}{c}{\multirow{2}{*}{\def\arraystretch{0.7}\begin{tabular}[c]{@{}c@{}}Time\\ per\\ iteration (s)\end{tabular}}} & \multicolumn{1}{c}{\multirow{2}{*}{\def\arraystretch{0.7}\begin{tabular}[c]{@{}c@{}}$\Delta$G\\ (kcal/mol)\end{tabular}}} & \multicolumn{1}{c}{\multirow{2}{*}{\def\arraystretch{0.7}\begin{tabular}[c]{@{}c@{}}Max\\ memory\\ usage (GB)\end{tabular}}} \\ \cline{2-7}
 & \multicolumn{1}{c}{LHS} & \multicolumn{1}{c}{RHS} & \multicolumn{1}{c}{\def\arraystretch{0.7}\begin{tabular}[c]{@{}c@{}}Precond-\\ itioning\end{tabular}} & \multicolumn{1}{c}{GMRES} & \multicolumn{1}{c}{\def\arraystretch{0.7}\begin{tabular}[c]{@{}c@{}}Calc.\\ of $\Delta$G\end{tabular}} & \multicolumn{1}{c}{Total} & \multicolumn{1}{c}{} & \multicolumn{1}{c}{} & \multicolumn{1}{c}{} & \multicolumn{1}{c}{} \\ \hline
\multicolumn{11}{c}{1bpi} \\ \hline
2 & 12.2 & 3.4 & 0.0 & 51.7 & 0.2 & 67.5 & 49 & 1.1 & -371.707 & 1.028 \\
10 & 18.8 & 4.7 & 0.3 & 47.3 & 0.3 & 71.4 & 18 & 2.6 & -371.420 & 1.141 \\
12 & 18.8 & 4.7 & 0.3 & 45.2 & 0.2 & 69.3 & 17 & 2.7 & -371.417 & 1.147 \\
20 & 15.9 & 4.6 & 5.4 & 59.7 & 0.2 & 85.8 & 11 & 5.4 & -371.722 & 1.149 \\
23 & 15.7 & 4.5 & 4.0 & 54.9 & 0.2 & 79.4 & 16 & 3.4 & -371.725 & 1.146 \\ \hline
\multicolumn{11}{c}{1lyz} \\ \hline
2 & 17.0 & 3.5 & 0.0 & 94.4 & 0.4 & 115.4 & 55 & 1.7 & -576.731 & 1.385 \\
10 & 27.5 & 5.0 & 0.8 & 81.4 & 0.4 & 115.1 & 19 & 4.3 & -576.071 & 1.595 \\
12 & 27.7 & 5.1 & 0.8 & 78.4 & 0.4 & 112.3 & 18 & 4.4 & -576.073 & 1.563 \\
20 & 24.4 & 5.1 & 6.9 & 113.4 & 0.4 & 150.2 & 13 & 8.7 & -576.761 & 1.601 \\
23 & 24.3 & 5.0 & 4.5 & 102.4 & 0.4 & 136.6 & 18 & 5.7 & -576.760 & 1.590 \\ \hline
\multicolumn{11}{c}{1a7m} \\ \hline
2 & 22.3 & 3.8 & 0.0 & 162.2 & 0.6 & 189.0 & 78 & 2.1 & -640.411 & 1.686 \\
10 & 38.8 & 5.6 & 1.3 & 104.6 & 0.6 & 150.8 & 20 & 5.2 & -640.176 & 2.008 \\
12 & 38.1 & 5.6 & 1.2 & 105.8 & 0.6 & 151.3 & 20 & 5.3 & -640.173 & 1.982 \\
20 & 34.4 & 5.7 & 8.8 & 158.5 & 0.6 & 208.0 & 15 & 10.6 & -640.416 & 1.995 \\
23 & 34.5 & 5.6 & 5.4 & 127.9 & 0.6 & 174.1 & 19 & 6.7 & -640.415 & 1.981 \\ \hline
\multicolumn{11}{c}{1x1z} \\ \hline
2 & 35.0 & 4.7 & 0.0 & 250.6 & 2.0 & 292.3 & 66 & 3.8 & -1669.642 & 2.779 \\
10 & 60.2 & 6.5 & 2.4 & 204.3 & 1.9 & 275.3 & 21 & 9.7 & -1667.484 & 3.288 \\
12 & 60.4 & 6.6 & 2.4 & 203.4 & 1.9 & 274.8 & 21 & 9.7 & -1667.450 & 3.289 \\
20 & 56.5 & 6.7 & 14.0 & 309.8 & 1.9 & 389.0 & 16 & 19.4 & -1669.467 & 3.286 \\
23 & 56.5 & 6.7 & 8.2 & 257.8 & 1.9 & 331.2 & 20 & 12.9 & -1669.473 & 3.267 \\ \hline
\multicolumn{11}{c}{4lgp} \\ \hline
2 & 63.6 & 5.6 & 0.1 & 600.3 & 6.8 & 676.4 & 85 & 7.1 & -2203.135 & 5.095 \\
10 & 115.5 & 9.2 & 8.9 & 365.5 & 6.7 & 505.9 & 20 & 18.3 & -2199.391 & 6.208 \\
12 & 114.3 & 9.2 & 8.6 & 394.4 & 6.7 & 533.1 & 22 & 17.9 & -2199.390 & 6.212 \\
20 & 110.7 & 9.1 & 28.8 & 690.4 & 6.7 & 845.7 & 19 & 36.3 & -2203.269 & 6.172 \\
23 & 110.2 & 9.2 & 17.2 & 503.1 & 6.7 & 646.3 & 21 & 24.0 & -2203.270 & 6.145 \\ \hline
\multicolumn{11}{c}{1igt} \\ \hline
2 & 117.7 & 8.6 & 0.1 & 1610.5 & 22.6 & 1759.5 & 126 & 12.8 & -4433.626 & 9.527 \\
10 & 218.7 & 15.2 & 13.5 & 743.6 & 22.7 & 1013.8 & 21 & 35.4 & -4428.218 & 11.558 \\
12 & 219.4 & 15.5 & 13.2 & 768.7 & 22.3 & 1039.2 & 23 & 33.4 & -4428.214 & 11.589 \\
20 & 210.5 & 15.0 & 46.9 & 1431.7 & 22.3 & 1726.3 & 22 & 65.1 & -4433.784 & 11.506 \\
23 & 214.5 & 15.2 & 25.1 & 1006.8 & 22.4 & 1284.0 & 24 & 41.9 & -4433.781 & 11.482 \\
\end{tabular}
}
\caption{Performance of the computational benchmarks.}
\label{tab:results}
\end{table}

In order to provide a detailed analysis of the formulations and preconditioners, we divide the simulation pipeline of the BEM into several phases, as reported in Table~\ref{tab:results}. First, the system matrix needs to be assembled, which is denoted as LHS (left-hand side). This operation only considers the singular parts of the matrix since all other interactions are calculated in each GMRES iteration with the FMM. Second, the right-hand-side (RHS) vector is constructed. This phase involves the calculation of the electrostatic potential and normal electric field from the atoms on the surface of the molecule, which was performed with a single FMM call. Third, the preconditioner is assembled. This involves the creation of the preconditioner matrix and, in the case of operator preconditioning, the sparse LU factorization of the mass matrices. Fourth, the GMRES algorithm solves the linear system, where the matrix-vector multiplications in each iteration are performed with the FMM. Finally, the solvation energy ($\Delta G$) calculation (Eq.~\eqref{eq:dG}), which includes the computation of the reaction potential (Eq.~\eqref{eq:phireac_bie}).

{\sffamily \small System assembly}\\

The assembly of the system matrix only involves the calculation of singular integrals, which are precomputed for the FMM. All other integrals are handled by the FMM in each GMRES iteration. The assembly time is expected to be directly related to the number of boundary integral operators in the specific formulation, which is confirmed by the column LHS in Table~\ref{tab:results}. The direct formulation only needs half the number of dense boundary integral operators as compared to the combined formulations. However, there is a discrepancy as the timings of the direct formulation are between one half and two thirds of the assembly time for the combined formulations. This is explained by two reasons. Firstly, the adjoint double-layer operator is not assembled since the transpose of the double-layer operator can be used. Secondly, the hypersingular operator that is present in the combined formulations, but not in the direct formulation, is more expensive to assemble than the single-layer and double-layer operators due to gradients in the operator kernel. Finally, the small reduction in assembly time of the first-kind formulation as compared to the second-kind formulation is because no identity operator is present.

{\sffamily \small Right-hand side calculation}\\

The assembly of the right-hand-side vector involves the calculation of the electrostatic potential and field on the molecular surface (see Eq.~\eqref{eq:sources}), and can be handled with the FMM. The direct formulation only needs the electrostatic potential while the combined field formulations also need the normal electric field from the sources. Hence, the calculation of the RHS is expected to be twice as long for the combined formulations, except for the overhead in setting up the FMM. This is confirmed by the third column in Table~\ref{tab:results}.

{\sffamily \small Preconditioning and matrix algebra}\\

In this study we used three preconditioning strategies: block-diagonal, (scaled) mass-matrix, and Calderón, which are all fast to create (fourth column in Table~\ref{tab:results}). In particular, the timings of the block-diagonal preconditioners are close to zero because all its information can be extracted from the previously assembled system matrix. The mass-matrix preconditioners need to calculate a sparse LU factorization of the discrete identity operator ($M$ in Eqs.~\eqref{eq:prec:mass} and \eqref{eq:prec:mass:scaled}). This becomes more expensive for larger molecules but the overhead is very small compared to the other phases of the software pipeline. The Calderón preconditioner uses the square of the system matrix, which is already available, and only the mass matrices that couple the preconditioner with the system matrix need to be assembled. Furthermore, the time to apply the preconditioner to the right-hand-side vector is also considered here. In the case of block-diagonal and mass-matrix, these are sparse matrix-vector products of the precomputed LU factorization, which are very fast. In the case of Calderón, this is done by a single FMM-accelerated dense matrix vector product, yielding a somewhat larger preconditioning time.

At each iteration, the GMRES performs a matrix-vector multiplication of the preconditioned matrix, accelerated by the FMM. Hence, the time per GMRES iteration strongly depends on the number of operators present in the system matrix and the preconditioner, as confirmed by the results in Table~\ref{tab:results}. The direct formulation has half the number of boundary integral operators compared to the combined field formulations. The matrix-vector multiplication is more than twice as fast because the combined field formulations use the hypersingular operator that requires a larger number of FMM evaluations than the single-layer and double-layer operators present in the direct formulation. 

Within Bempp-cl, a matrix-matrix-vector multiplication is implemented as two consecutive matrix-vector products. Then, the Calderón preconditioned system (Case 20) involves an extra (dense) matrix-vector product at each GMRES iteration, which explains why the time per iteration is double compared to second-kind formulations. In comparison, the preconditioner in Case 23 has half the number of operators, leading to a decrease in the total time per iteration. 

{\sffamily \small GMRES convergence}\\

\begin{figure}[!ht]
    \centering
    \includegraphics[width=0.7\columnwidth,keepaspectratio=true]{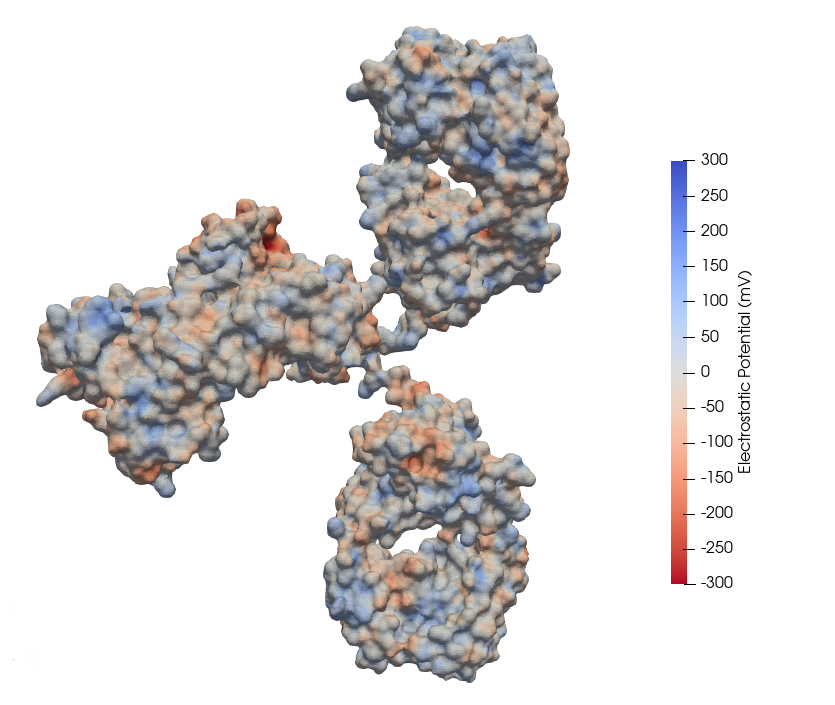}
    \caption{\label{fig:1igt} The electrostatic potential on the surface of the 1igt molecule.}
\end{figure}

The GMRES algorithm solves the linear system resulting from the different formulations, to obtain the electrostatic potential and normal electric field on the molecular surface. As an example, Figure~\ref{fig:1igt} shows the electrostatic potential on the solvent-excluded surface of immunoglobulin G (PDB code 1igt), retrieved directly from the output of the GMRES. The benchmarks show that the FMM-accelerated linear solver is the most expensive part of the algorithm, consuming up to 90\% of the total time to solution (see columns 5 and 7 in Table~\ref{tab:results}). Hence, the design of efficient preconditioners to reduce the number of GMRES iterations has a significant impact on the overall efficiency of the implicit-solvent model. Analyzing column~8 in Table~\ref{tab:results}, it is clear that the direct formulation is the worst conditioned, which is expected because it has no eigenvalue clustering. Indeed, Case 2 requires a larger number of GMRES iterations compared to the combined field formulations, to the point that, for example, 1igt has six times more iterations. Even worse, the convergence deteriorates quickly with the problem size. Differently, the second-kind formulations are well-conditioned by design and converge in a few GMRES iterations. When comparing the smallest with the largest molecule, the number of degrees of freedom increases with a factor of nine but the iteration count only from 18 to 21 for the scaled mass-matrix preconditioned Juffer formulation (Case 10). Furthermore, Calderón preconditioning of the first-kind formulation is very effective, with only 11 iterations needed for 1bpi (Case 20), the lowest number of any preconditioned formulation. However, the convergence deteriorates moderately with problem size, reaching the same level as the second-kind formulations for the 1igt molecule, thus loosing its competitive advantage. This is specially worrying since the time per iteration is twice as long. Finally, the internal Calderón preconditioner (Case 23) shows the design challenges of preconditioning strategies. That is, it requires less but slower iterations than the mass-preconditioned formulations (Cases 10 and 12), and more but faster iterations than the full Calderón preconditioner (Case 23).

{\sffamily \small Solvation energy}\\

The sixth column in Table \ref{tab:results} shows that the calculation of the solvation energy ($\Delta G$) is small compared to the other subroutines. This step involves computing the reaction potential ($\phi_\text{reac}$) at the location of the atoms from the potential at the surface (see Eq.~\eqref{eq:phireac_bie}), and then integrating this value with the Dirac-delta charge distribution (see Eq.~\eqref{eq:dG}). The benchmarks confirm that the time taken to compute $\Delta G$ is agnostic to the choice of preconditioned formulation, as expected.

{\sffamily \small Total computation time}\\

The results in Table~\ref{tab:results} show that the total time is a balance between a reduced number of GMRES iterations at the expense of longer times per iteration. Then, no preconditioned formulation is expected to be the best (in terms of calculation time) for all types of molecules. Comparing the total time of the different models, two general recommendations can be made.

First, the direct formulation with block-diagonal preconditioning (Case~2) is the quickest algorithm for smaller molecules. As expected, the assembly of the system and each GMRES iteration are very fast since only the identity, single-layer and double-layer operators are present. Furthermore, the block-diagonal preconditioner has almost no computational overhead. Since small molecules have a low number of mesh elements, the GMRES converges quickly for any reasonable formulation anyway. Hence, the fast arithmetic of the direct formulation makes it the best choice for smaller molecules.

Second, the second-kind formulations with scaled mass-matrix preconditioning are the quickest algorithms at larger molecules, with the Juffer formulation (Case~10) slightly faster than the Lu formulation (Case~12). The mathematical analysis of the original second-kind formulations predicts a clustering of eigenvalues around accumulation points. Furthermore, the scaled mass-matrix preconditioning yields a single accumulation point of the preconditioned system. Hence, the preconditioned formulations are expected to be well-conditioned for a wide range of problem sizes. Even though the assembly and each GMRES iteration are more expensive than those of the direct formulation, the GMRES converges six times faster. This confirms the importance of designing well-conditioned formulations and effective preconditioners to solve the implicit-solvent model at larger molecules.

\subsection*{\sffamily \large Fast Calderón preconditioning}

\begin{table}
\centering
\resizebox{\textwidth}{!}{%
\def\arraystretch{0.8}
\begin{tabular}{crrrrrrrrrr}
\hline
\multirow{2}{*}{\def\arraystretch{0.7}\begin{tabular}[c]{@{}c@{}}Case\\ N°\end{tabular}} & \multicolumn{6}{c}{Timings (s)} & \multicolumn{1}{c}{\multirow{2}{*}{\def\arraystretch{0.7}\begin{tabular}[c]{@{}c@{}}GMRES\\ iterations\end{tabular}}} & \multicolumn{1}{c}{\multirow{2}{*}{\def\arraystretch{0.7}\begin{tabular}[c]{@{}c@{}}Time\\ per\\ iteration (s)\end{tabular}}} & \multicolumn{1}{c}{\multirow{2}{*}{\def\arraystretch{0.7}\begin{tabular}[c]{@{}c@{}}$\Delta$G\\ (kcal/mol)\end{tabular}}} & \multicolumn{1}{c}{\multirow{2}{*}{\def\arraystretch{0.7}\begin{tabular}[c]{@{}c@{}}Max\\ memory\\ usage (GB)\end{tabular}}} \\ \cline{2-7}
 & \multicolumn{1}{c}{LHS} & \multicolumn{1}{c}{RHS} & \multicolumn{1}{c}{\def\arraystretch{0.7}\begin{tabular}[c]{@{}c@{}}Precond-\\ itioning\end{tabular}} & \multicolumn{1}{c}{GMRES} & \multicolumn{1}{c}{\def\arraystretch{0.7}\begin{tabular}[c]{@{}c@{}}Calc.\\ of $\Delta$G\end{tabular}} & \multicolumn{1}{c}{Total} & \multicolumn{1}{c}{} & \multicolumn{1}{c}{} & \multicolumn{1}{c}{} & \multicolumn{1}{c}{} \\ \hline
\multicolumn{11}{c}{1bpi} \\ \hline
20b & 15.6 & 4.6 & 12.0 & 40.9 & 0.3 & 73.3 & 12 & 3.4 & -371.722 & 1.229 \\
23b & 15.7 & 4.5 & 8.8 & 45.5 & 0.3 & 74.8 & 16 & 2.8 & -371.726 & 1.163 \\ \hline
\multicolumn{11}{c}{1lyz} \\ \hline
20b & 24.8 & 5.1 & 16.9 & 74.8 & 0.4 & 122.0 & 14 & 5.3 & -576.759 & 1.697 \\
23b & 24.8 & 5.0 & 11.4 & 88.6 & 0.4 & 130.2 & 19 & 4.7 & -576.761 & 1.618 \\ \hline
\multicolumn{11}{c}{1a7m} \\ \hline
20b & 34.0 & 5.6 & 22.3 & 97.9 & 0.6 & 160.4 & 15 & 6.5 & -640.416 & 2.119 \\
23b & 34.4 & 5.6 & 13.5 & 105.0 & 0.6 & 159.1 & 19 & 5.5 & -640.416 & 2.019 \\ \hline
\multicolumn{11}{c}{1x1z} \\ \hline
20b & 55.3 & 6.6 & 34.0 & 198.3 & 1.9 & 296.1 & 16 & 12.4 & -1669.475 & 3.521 \\
23b & 55.9 & 6.6 & 20.1 & 221.0 & 1.9 & 305.5 & 21 & 10.5 & -1669.473 & 3.342 \\ \hline
\multicolumn{11}{c}{4lgp} \\ \hline
20b & 110.6 & 9.1 & 67.0 & 464.3 & 6.7 & 657.8 & 21 & 22.1 & -2203.271 & 6.838 \\
23b & 112.9 & 9.2 & 38.6 & 452.6 & 6.7 & 620.1 & 23 & 19.7 & -2203.270 & 6.445 \\ \hline
\multicolumn{11}{c}{1igt} \\ \hline
20b & 211.5 & 14.9 & 121.3 & 909.9 & 22.4 & 1280.0 & 23 & 39.6 & -4433.784 & 12.776 \\
23b & 212.7 & 14.9 & 65.6 & 889.3 & 22.2 & 1204.7 & 26 & 34.2 & -4433.781 & 11.971
\end{tabular}
}
\caption{Performance of the computational benchmarks for fast Calderón preconditioning.}
\label{tab:results-fast-calderon}
\end{table}

The computational benchmarks presented in Table~\ref{tab:results} show that the Calderón preconditioned first-kind formulation is too slow compared to the direct formulation and too ill-conditioned to outperform the second-kind formulations. The expensive preconditioning strategy of squaring the matrix doubles the time per iteration, and the number of GMRES iterations scales unfavorably with problem size. Hence, it is the slowest in almost any benchmark.

An improvement of the Calderón preconditioner is to consider different numerical parameters for the preconditioner than for the system matrix~\cite{escapil2019fast}. In particular, we can use relaxed quadrature and FMM parameters on the preconditioner only, without sacrificing overall accuracy of the solution. In this case, we reduced the quadrature order from 3 to 1 for the regular operators and from 4 to 3 for the singular operators, and reduced the FMM expansion order from 3 to 2. Table~\ref{tab:results-fast-calderon} shows results for this strategy, tagged as Cases 20b (full Calderón) and 23b (scaled-mass internal Calderón).

This fast Calderón preconditioning strategy leads to faster matrix-vector products for the preconditioner, yielding a lower time per iteration, but at the expense of longer assembly time since the system matrix cannot be reused anymore. The benchmarks confirm that the assembly of the preconditioner takes much longer than before, even though not as slow as creating the system matrix. The advantage is that each GMRES iteration only takes roughly 60\% of the standard Calderón preconditioner. Notice that we chose the internal-Calderón preconditioner over the external-Calderón variant since the interior operators for the Laplace equation are quicker with the FMM compared to the linearized Poisson-Boltzmann equation in the exterior medium (also known as the Yukawa kernel). Furthermore, the preconditioner is less accurate, which leads to slower convergence of GMRES, but the increase in iteration count is small. Even though the fast Calderón preconditioner does not outperform the other formulations in overall timings, for small and intermediate molecules, the GMRES time is better than any other formulation. This is beneficial when the system needs to be solved for multiple right-hand-side vectors, for example, in pKa calculations~\cite{potter1994small,warwicker2011pka,hu2018accurate,chen2021computing}.

Notice that one can take this strategy further and consider, for example, only the near field interactions of the FMM or single-precision arithmetic for the preconditioner. However, these changes require invasive adjustments to the algorithm and software. This is outside the scope of this study because we aimed to solve the implicit-solvent model with readily available implementations and easy coding strategies. In any case, the benchmarks show a deterioration of the Calderón preconditioning for large molecules. The number of iterations for the largest molecule (1igt) is roughly the same as for the Juffer and Lu formulations. Hence, improving the time per GMRES iteration for the Calderón preconditioner will not yield better overall timing than the scaled mass-matrix preconditioned Juffer formulation at large-scale simulations.

\subsection*{\sffamily \large Notes on implementation}

All formulations were implemented using the Bempp-cl library~\cite{betcke2021bempp} on the \texttt{bem\_electrostatics} code\footnote{\url{https://github.com/bem4solvation/bem_electrostatics}}. This code is easily downloaded and installed via \texttt{conda} and \texttt{pip}, and can be called from a Jupyter notebook with a single command, for different choices of formulations and preconditioners. All the required instructions to reproduce the results are available in a Github repository\footnote{\url{https://github.com/bem4solvation/formulations_paper}}.

The level of abstraction of the boundary integral operators through a Python interface offered by Bempp-cl makes it easy to prototype and test the different combinations of formulations and preconditioners. This way, we could do a thorough search of cases with minor implementation overhead, on a platform that allows us to reach moderately large scales in reasonable time. This is a key aspect of this work, especially as we solve bigger problems, where the computational efficiency of the method is more critical. This infrastructure becomes useful to explore the different models before simulating at extremely large scales~\cite{martinez2019free}.


 





\section*{\sffamily \Large CONCLUSIONS}
This paper presents a generalized boundary integral formulation of the Poisson-Boltzmann equation for an implicit-solvent model, and proposes a methodology to construct efficient preconditioners for each case. Well-known boundary integral formulations of this problem, such as Juffer's second-kind equation, are special cases of our generalized approach, and it allows us to easily create and test new expressions. We studied the performance of a set of formulations and preconditioners for a small test case with arginine, focusing on the time to solution and matrix conditioning, through the eigenvalue spectrum. We saw that a simple scaling of the mass-matrix preconditioner allowed us to manipulate and cluster the eigenvalues, improving the conditioning of the system. Then, we chose the best performing combinations of formulation and preconditioner and assessed their scaling with the solute size, up to $\sim$220,000 degrees of freedom. We found that the optimal choice is problem-size dependent, and is a balance between the number of GMRES iterations and the speed of the FMM-accelerated matrix-vector product in each iteration. In particular, the direct method was the fastest for small molecules, whereas second-kind formulations were most efficient as the molecule size grew. We also described a first-kind expression with a fast Calderón preconditioner that even though it was slower than second-kind formulations, it presented promising properties such as short GMRES times. We expect this study to be useful for researchers to find the best strategy to compute the electrostatic potential and its contribution to the solvation energy in molecular modelling, specially as they aim at large-scale simulations.

As future work, we plan to explore more formulations, and use other comparison criteria. For example, there is evidence that first-kind formulations can be more accurate than second-kind equations. Then, smaller mesh densities solved with a first-kind form would yield equivalent precision to a second-kind counterpart, and may result in a faster method. We will also continue to study the fast Calderón preconditioner, looking to increase the speed by using more invasive adjustments of the software.


\subsection*{\sffamily \large ACKNOWLEDGMENTS}
SDS and CDC acknowledge the support from Universidad Técnica Federico Santa María through project PI-LIR-2020-10. 


\clearpage


\bibliography{bibtexrefs}   

\end{document}